\documentclass[4paper,11pt]{amsart}

\usepackage{style}

\title[Induced birational transformations on O'Grady's sixfolds]{Induced birational transformations on O'Grady's sixfolds}

\author{Annalisa Grossi}
\address{Annalisa Grossi, Fakult\"{a}t f\"{u}r Mathematik, Technische  Universit\"{a}t Chemnitz, Reichenhainer Strasse 39, 09126 Chemnitz, Germany}
\email{annalisa.grossi@math.tu-chemnitz.de}

\date{\today}

\makeatletter
\@namedef{subjclassname@2020}{%
 \textup{2020} Mathematics Subject Classification}
\makeatother


\thanks{}

\usepackage{framed}

\begin{document}

\maketitle

\begin{abstract}
	We introduce the notion of induced birational transformations of irreducible holomorphic symplectic sixfolds of the sporadic deformation type discovered by O'Grady. We give a criterion to determine when a manifold of \(\OG_6\) type is birational to \(\widetilde{\K}_v(\A, \theta)\), a moduli space of sheaves on an abelian surface. Then we determine when a birational transformation of \(\widetilde{\K}_v(\A, \theta)\) is induced by an automorphism of \(\A\). Referring to the Mongardi--Rapagnetta--Sacc\'{a} birational model of manifolds of \(\OG_6\) type, we give a result to determine when a birational transformation is induced at the quotient. We give an application of these criteria in the nonsymplectic case.
\end{abstract}

\section{Introduction}

This paper deals with induced birational transformations of irreducible holomorphic symplectic sixfolds of O'Grady's deformation type. Throughout the paper we will refer to sixfolds of O'Grady's type as manifolds of \(\OG_6\) type. Irreducible holomorphic symplectic manifolds arise from symplectic surfaces, and in many cases are constructed as moduli spaces of sheaves on them. When we consider automorphisms, or more in general, birational transformations of irreducible holomorphic symplectic manifolds is then natural to ask weather they are induced by an automorphism of the K3 surface involved in the construction of the moduli space.

Hilbert schemes of \(n\) points on a K3 surface \cite{beauville1983varietes} allows Beauville \cite{beauville1983some} to extend several classical results of Nikulin and Boissiere \cite{boissiere2012automorphismes} to introduce the notion of \textit{natural} automorphism simply by taking a nontrivial automorphism of the K3 and considering the induced action on its Hilbert scheme. 
A generalization of the notion of natural automorphisms for moduli spaces of sheaves on symplectic surfaces appeared for the first time in a work of Ohashi--Wandel \cite{ohashi2013non}, which is inspired by a construction due to Oguiso--Schr\"{o}er  \cite{oguiso2011enriques}.
This notion was adapted to more general cases by Mongardi--Wandel \cite{MW15}, using recent developments in stability condition theory due to Bridgeland \cite{bridgeland2008stability}, Bayer-Macrì \cite{bayer2014mmp},\cite{bayer2014projectivity} and Yoshioka \cite{yoshioka2016bridgeland}.
Moreover, in \cite{MW15}, Mongardi--Wandel conjecture the possibility to extend the notion of induced automorphisms also in the case of manifolds of O'Grady's deformation type.

In this paper we give a notion of induced automorphisms and, more in general, of induced birational transformations of irreducible holomorphic symplectic manifolds of \(\OG_6\) type, considering two birational models of them.  The first birational model is the resolution of the Albanese fiber of a moduli space of sheaves on an abelian surface. The second model due to Mongardi--Rapagnetta--Sacc\`{a} is the resolution of the quotient of the Hilbert scheme of \(3\) points on a K3 surface by a birationl symplectic involution. We introduce the notions of induced birational transformation and birational transformation induced at the quotient. The first notion refers to the first birational model and it means essentially that a birational transformation of the manifold of \(\OG_6\) type comes from an automorphism of the abelian surface. In the second case the notion refers to the second birational model and it means that a birational transformation of the manifold of \(\OG_6\) type can be lifted to a birational transformation of the Hilbert scheme involved in the construction.

\subsection{Contents of the paper}
In \autoref{preliminaries} we introduce basic tools of lattice theory for irreducible holomorphic symplectic manifolds, and we recall the construction of O'Grady's sixfolds, due to the contribution of many authors: O'Grady \cite{o2000new}, Kaledin--Lehn--Sorger \cite{Kal_Lehn_Sorger_singular}, Lehn--Sorger \cite{Lehn_Sorger} and Perego--Rapagnetta \cite{perego2013deformation}. O'Grady \cite{o2000new} introduces the sporadic example in dimension six as the symplectic resolution of the fibre of an isotrivial fibration defined on a moduli space of sheaves on an abelian surface \(\A\), with respect to a non-primitive Mukai vector \(v\) and a \(v\)-generic polarization \(\theta\), and denotes this fibre by \(\widetilde{\K}_v(\A, \theta)\). Later Perego--Rapagnetta \cite{perego2013deformation} generalize this construction. They introduce the notion of OLS-triple and they find that also with more general assumptions on \(v\), \(\A\) and \(\theta\), the O'Grady's construction holds true and the manifold \(\widetilde{\K}_v(\A, \theta)\) is deformation equivalent to the O'Grady's six dimensional example. 

In \autoref{Induced automorphisms groups} we give a lattice-theoretic criterion to determine when a manifold of \(\OG_6\) type is birational to \(\widetilde{\K}_v(\A,\theta)\). It is a general fact that the second integral cohomology group of irreducible holomorphic symplectic manifolds is endowed with a lattice structure. We introduce the notion of \textit{numerical moduli space} (see \autoref{nummodspace}) for a manifold \(X\) of \(\OG_6\) type, which concerns in some conditions on the second integral cohomology lattice of \(X\). A marked pair \((X,\eta)\) of \(\OG_6\) type is a pair where \(X\) is a manifold of \(\OG_6\) type, and \(\eta \colon \HH^{2}(X,\mathbb{Z}) \to \bL\) is a fixed isometry of lattices. More precisely, we give the following characterization.

\begin{theorem}\label{X n.m.s allora birazionale alla fibra della mappa di albanese}
If \((X, \eta)\) is a marked pair of \(\OG_6\) type, then there exists an abelian surface \(\A\), a non-primitive Mukai vector \(v=2w\), and a \(v\)-generic polarization \(\theta\) on \(\A\) such that \(X\) is birational to \(\widetilde{\K}_v(\A, \theta)\) if and only if \(X\) is a numerical moduli space.
\end{theorem}

See \autoref{proof of Thm 1.1} for the proof.
Then we give the definition of \textit{numerically induced group of birational transformations} (see \autoref{numerically induced}) and we prove the following theorem to determine when a birational transformation of a manifold of \(\OG_6\) type is induced by an automorphism of the abelian surface. We denote by \(\Bir(X)\) the group of birational transformations of \(X\).

\begin{theorem}\label{num ind impliesd ind}
Let \((X, \eta)\) be a smooth marked pair of \(\OG_6\) type. Let \(\G \subset \Bir(X)\) be a finite subgroup. If \(\G\) is a numerically induced group of birational transformations then there exists an abelian surface \(\A\) with \(\G \subset \Aut(\A) \times \A^{\vee}[2]\), a \(\G\)-invariant Mukai vector \(v\) and a \(v\)-generic polarization \(\theta\) on \(\A\) such that \(X\) is birational to \(\widetilde{\K}_v(\A,\theta)\) and \(\G\subset \Bir(X)\) is an induced group of birational transformations.
\end{theorem}
See \autoref{section_4.1} for the proof.

Mongardi--Rapagnetta--Sacc\`{a} \cite{MRSHodge} prove that there exist manifolds of \(\OG_6\) type that admit a birational model obtained as a quotient of a manifold of K3\(^{[3]}\) type by a birational symplectic involution. In \autoref{Automorphisms induced at the quotient} we recall the construction of Mongardi--Rapagnetta--Sacc\`{a} and we give a lattice-theoretic criterion to determine when a birational transformation of a manifold of \(\OG_6\) type lifts to a birational transformation of the manifold of K3\(^{[3]}\) type involved in the construction. In such a case we call the birational transformation induced at the quotient (see \autoref{def aut ind at the quotient}). 

\begin{theorem}\label{ind at the quot}
Let \(X\) be a manifold of \(\OG_6\) type which is a numerical moduli space and let \(\varphi \in \Bir(X)\) be a birational transformation of \(X\). If there exists a class \(E \in \NS(X)\) of square \(-2\) and divisibility \(2\) which is fixed by the induced action of \(\varphi\) in cohomology, then \(\varphi\) is induced at the quotient.
\end{theorem}

See \autoref{proof thm 1.3} for the proof. In \autoref{A sufficient cond to have a regular morphism of K3} we prove \autoref{thm induced at the quotient} which state a sufficient condition to extend the birational transformation to an automorphism of the manifold of K3\(^{[3]}\) type. The condition that we state is geometric and it is related to the action of the induced action on the singular locus of the \(2 \colon 1\) cover of the singular moduli space \(\K_v(\A, \theta)\). 

Finally in \autoref{Applications} we apply our techniques to nonsymplectic automorphisms (automorphisms that do not preserve the symplectic form) of manifolds of \(\OG_6\) type. Using a classification  of nonsymplectic automorphisms of prime order of manifolds of \(\OG_6\) type contained in \cite{Grossi:non-sp.aut.OG3} we prove the following theorem.

\begin{theorem}\label{theorem applications}
Let \((X, \eta)\) be a marked pair of \(\OG_6\) type. Assume that \(X\) is a numerical moduli space and let \(\varphi \in \Aut(X)\) be a nonsymplectic automorphism of \(X\) of prime order. Then \(\varphi\) is induced and induced at the quotient in cases that are listed in \autoref{tab:L}.
\end{theorem}
A complete proof of it is given in \autoref{proof of thm 1.4}. In \autoref{induced implies induced at the quotient} we show that if a birational transformation is induced then it is induced at the quotient, and in \autoref{tab:L} we denote by \(\clubsuit\) the involutions that are induced at the quotient but not induced. By \autoref{S_G(X) should be even in tha p=2 case} nonsymplectic automorphisms of prime order that can be induced at the quotient but not induced are among involutions. 
\begin{center}
\small{
	\begin{longtable}{llllllll}
		
	\caption{Induced and induced at the quotient groups of nonsymplectic automorphisms of prime order on manifolds of \(\OG6\) type.}
		\label{tab:L} \\
		
		\toprule
	No. & \(|\G|\)  & \(\bL_G\) & \(\bL^G\) & induced & {\small{ind. at the quot.}} &  \\

	\midrule
	\endfirsthead
	
	\multicolumn{5}{c}%
	{\tablename\ \thetable{}, follows from previous page} \\
	\midrule
		No. & \(|\G|\) ) & \(\bL_G\) & \(\bL^G\) & induced & ind. at the quot. &  \\
	\midrule
	\endhead
	
	\midrule
	\multicolumn{7}{c}{Continues on next page} \\
	\endfoot
	
	\bottomrule
	\endlastfoot
	
	\(1\) & \(2\)  & \(\bU^{\oplus 2} \oplus [ -2 ]^{\oplus 3 }\)& \([2]\) & no & no &\\
	\(2\) & \(2\)   & \(\bU \oplus [2] \oplus [-2] ^{\oplus 3 }\) & \([2] \oplus [-2]\) & no & no & \\
	\(3\) & \(2\) & \(\bU^{\oplus 2} \oplus [-2 ] ^{\oplus 2 }\) & \(\bU\) & no & no  &\\
	\(4\) & \(2\)  & \(\bU^{\oplus 2} \oplus [-2 ] ^{\oplus 2 }\) & \([ 2 ] \oplus [-2]\) & no & no  &\\
	\(5\) & \(2\)  & \(\bU^{\oplus 2} \oplus [-2 ] ^{\oplus 2 }\) & \(\bU(2)\) & no & no &\\
	\(6\) & \(2\)   & \([2]^{\oplus2} \oplus [-2]^{\oplus3}\) & \([2] \oplus[-2]^{\oplus2}\) & no  & no &\\
	\(7\) & \(2\)  & \(\bU \oplus [-2] ^{\oplus 2 } \oplus [2]\) & \(\bU \oplus [-2]\) & no & no & \\
	\(8\) & \(2\)  & \(\bU \oplus [-2] ^{\oplus 2 } \oplus [2]\) & \([2] \oplus [-2]^{\oplus 2}\) & no & no & \\
	\(9\)  & \(2\) & \(\bU^{\oplus 2} \oplus [-2]\) & \([-2] ^{\oplus 2} \oplus [2]\) & no & yes  & \(\clubsuit\) \\
	\(10\)  & \(2\)  & \(\bU^{\oplus 2} \oplus [-2]\) & \(\bU \oplus [-2]\) & no  & yes & \(\clubsuit\) \\
	\(11\) & \(2\)  & \([2] ^{\oplus 2} \oplus [-2] ^{\oplus 2}\) & \(\bU \oplus [-2] ^{\oplus 2}\) & no & no &\\
  \(12\) & \(2\)  & \([2]  \oplus[-2] \) & \(\bU \oplus [2] \oplus [-2] ^{\oplus 3}\) & yes & yes &\\
	\(13\) & \(2\) & \(\bU(2)^{\oplus 2}\)  & \(\bU(2) \oplus [-2]^{\oplus2}\) & no  & no &\\
	\(14\) & \(2\)   & \(\bU \oplus [2] \oplus [-2]\) & \([2] \oplus [-2]^{\oplus 3}\) & yes & yes &\\
	\(15\) & \(2\)& \(\bU \oplus [2] \oplus [-2]\) & \(\bU \oplus [-2] ^{\oplus 2}\) & yes & yes &\\
	\(16\) & \(2\)  & \(\bU \oplus \U(2)\) & \(\bU(2) \oplus [-2] ^{\oplus 2}\) & yes & yes &\\
	\(17\) & \(2\) & \(\bU \oplus \bU(2)\) & \(\bU \oplus [-2] ^{ \oplus 2}\) & yes & yes & \\
	\(18\) & \(2\)  & \(\bU ^{\oplus 2}\)& \(\bU \oplus [-2] ^{\oplus 2}\) & yes  & yes &\\
	\(19\) & \(2\)   & \([2]^{\oplus 2} \oplus [-2]\) & \([-2]^{\oplus 4} \oplus [2]\) & no  & yes & \(\clubsuit\) \\
	\(20\) & \(2\)  & \([2]^{\oplus 2} \oplus [-2]\) & \(\bU \oplus [-2]^{\oplus 3}\) & no  & yes & \(\clubsuit\) \\
	\(21\)& \(2\)  & \(\bU \oplus [2]\) & \(\bU \oplus [-2] ^{\oplus 3}\) & no & yes & \(\clubsuit\)\\
    \(22\) & \(2\) & \([2]^{ \oplus 2}\)& \(\bU \oplus[-2] ^{\oplus 4}\) & yes & yes &\\
			\(23\) & \(2\)& \([2]^{ \oplus 2}\)& \(\bU \oplus \bD_{4}(-1)\)& no & no & \\
			\(24\) & \(2\) & \([2]^{\oplus2}\) & \(\bU(2) \oplus \bD_{4}(-1)\) &  no & no & \\
			\midrule

			\(1\) & \(3\)  & \(\bU^{\oplus2} \oplus \bA_{2}(-1)\) & \([-2] \oplus [6]\) & no & no \\
			\(2\) & \(3\) & \(\bA_{2}\) & \(\bU \oplus \bA_{2}(-1) \oplus [-2]^{\oplus 2}\) & yes & yes &\\
			\midrule
			
			\(1\) & \(5\)  & \(\bU \oplus \h_{5}\) & \([-2] \oplus [-10] \oplus \bU\) & yes & yes & \\
			\midrule
			
			\(1\) & \(7\)  & \(\bU^{\oplus 2} \oplus \Bk_{7}\) & \([-2] \oplus [14]\) & no & no &\\
			
	\end{longtable}}
\end{center}

\textbf{Acknowledgments}. I would like to warmly thank
Samuel Boissi\`{e}re, Giovanni Mongardi, Claudio Onorati, Antonio Rapagnetta, Alessandra Sarti and Davide Cesare Veniani
for the many many useful hints and for fruitful discussions.


\section{Preliminaries}\label{preliminaries}
In this section we fix the notation and the conventions that we will use throughout the paper. Moreover we recall basic results for irreducible holomorphic symplectic manifolds of \(\OG_6\) type, and we present them in a more suited form to the purposes of this work.

In \autoref{lattices} we recall some definitions and results of lattice theory for irreducible holomorphic symplectic manifolds and, in particular, for manifolds of \(\OG_6\) type, and we collect some basic results about primitive embeddings of lattices. In \autoref{O'Grady's sixfold} we summarize the original construction of O'Grady's sixfolds \cite{o2000new} and the more general one due to Perego--Rapagnetta \cite{perego2013deformation}, which construct manifolds of \(\OG_6\) type as moduli spaces of sheaves on abelian surfaces. Moreover we recall how it is possible to endow the second integral cohomology of these moduli spaces with a pure weight--two Hodge structure. 
\subsection{Lattice theory for irreducible holomorphic symplectic manifolds}\label{lattices}
\subsubsection{Lattices}\label{Lattices_SUB}
A \textit{lattice} \(L\) is a free \(\mathbb{Z}\)-module of finite rank endowed with a nondegenerate symmetric bilinear form \[
    L \times L \rightarrow \mathbb{Z}, 
    \]
    \[
(e,f) \mapsto e\cdot f.
\]
We denote \(e \cdot e\) by \(e^{2}\).
The lattice \(L\) is \textit{even} if \(e^2 \in 2\mathbb{Z}\).
Every lattice has a determinant, a signature and a rank denoted by \(\det(L)\), \(\sign(L)\) and \(\rank(L)\) respectively.
The \textit{dual} of a lattice \(L\) is \(L^{\vee}=\Hom_{\mathbb{Z}}(L,\mathbb{Z})\), which admits the following equivalent description:
\[
L^{\vee}=\{x \in L \otimes \mathbb{Q} \ |\  x \cdot y \in \mathbb{Z} \mbox{ for each } y \in L\}.
\]
To every lattice \(L\) we can associate a finite group \(L^{\sharp}=L^{\vee}/L\), which is called \textit{discriminant group} of \(L\). The order of \(L^{\sharp}\) is  \(|{\det(L)}|\). The \textit{length} of the discriminant group \(L^{\sharp}\) is the minimal number of generators of it and it is denoted by \(l(L^{\sharp})\).
A lattice \(L\) is said to be \textit{unimodular} if \(L^{\sharp}=\{\id\} \) and \textit{p-elementary} if \(L^{\sharp}=(\mathbb{Z}/p\mathbb{Z})^{\oplus a}\) for a prime number \(p\) and \(a \in \mathbb{Z}_{\geq 0}\). In this case \(l(L^{\sharp})=a\). To avoid any ambiguity for \(a=0\) we include here the case of unimodular lattices, that are considered \(p\)-elementary for any \(p\).

The pairing on an even lattice \(L\) induces a pairing on \(L^{\sharp}\) with values in \(\mathbb{Q}/\mathbb{Z}\) and the associated \(\mathbb{Q}/2\mathbb{Z}\)-valued quadratic form is 
\[
q_{L}: L^{\sharp} \to \mathbb{Q}/2\mathbb{Z},
\]
\[
q_{L}(e+L)=(e,e) \ \ \mbox{mod } 2\mathbb{Z}
\]
that is called \textit{discriminant form of L}.
 There exists a natural homomorphism \(\OO(L) \to \OO(L^{\sharp})\). We denote by \(G^{\sharp}\) the image in \(\OO(L^{\sharp})\) of a subgroup of isometries \(G \subset \OO(L)\).
 The \textit{divisibility} of a vector \(v \in L\) is defined as \((v,L)=\mbox{gcd}\{v \cdot v' \ | \ v' \in L\}\) and \(v/(v,L)\) is an element of \(L^{\vee}\), hence of \(L^{\sharp}\). In this paper we will refer to the unique even unimodular indefinite lattice of rank two by \(\bU\). Moreover \(\bA_n\), \(\bD_n\) and \(\bE_n\) denote the negative definite ADE lattices. The notation \([m]\) with \(m \in \mathbb{Z}\) refers to a lattice of rank \(1\) generated by a vector of square \(m\). If \(L\) is a lattice, we denote by \(L(n)\) the lattice with the same structure as \(\mathbb{Z}\)-module, but with quadratic form multiplied by \(n\).
 
 Two lattices have the same \textit{genus} if they have the same signature and isomorphic quadratic forms. Due to Nikulin's criteria \cite[Theorem 3.6.2]{Nik79}, \cite[Theorem 1.14.2]{Nik79} for indefinite lattices, and due to the Smith--Minkowski--Siegel formula for definite lattices \cite{conway_low_dim}, we know that the lattices that appear in this paper are unique in their genus.
 
  If \(L\) is a lattice which is endowed with a pure weight-two Hodge structure and if we denote \(L \otimes_{\mathbb{Z}} \mathbb{C}\) by \(L_{\mathbb{C}}\) then
\[ L_{\mathbb{C}}=L_{\mathbb{C}}^{2,0} \oplus L_{\mathbb{C}}^{1,1} \oplus L_{\mathbb{C}}^{0,2}.
\]
 Moreover it holds that \(\overline{L_{\mathbb{C}}^{2,0}}=L_{\mathbb{C}}^{0,2}\), \(\overline{L_{\mathbb{C}}^{1,1}}=L_{\mathbb{C}}^{1,1}\) and \(\dim_{\mathbb{C}}(L_{\mathbb{C}}^{2,0})=1\). In this paper we use the notation 
\[
L^{1,1}=L_{\mathbb{C}}^{1,1} \cap L,
\]
to refer to the integral \((1,1)\)-part of the lattice \(L\).

\subsubsection{Primitive embeddings and isometries}\label{Primitive embeddings and isometries} An embedding of lattices \(M \subset L\) is called primitive if the group \(L/M\) is torsion-free. In this setting we denote the embedding by \(M \hookrightarrow L\) and we denote by \(N=M^{\perp}\) the \textit{orthogonal complement} of \(M\) in \(L\). To compute the possible primitive embeddings of \(M\) in \(L\) we refer to \cite[Proposition 1.5.1]{Nik79}. The definitions and the notations that we recall here refer to \cite[Section 2.2]{gro_on_ven2021symplectic}.

If there exists a primitive embedding of even lattices \(M \hookrightarrow L\) then there exists a subgroup, called a \textit{gluing subgroup} \(H \subset M^{\sharp}\) 
and a gluing isometry \(\gamma \colon H \to H' \subset N(-1)^{\sharp}\). If \(\Gamma\) denotes the gluing graph of \(\gamma\) in \(M^{\sharp} \oplus N^{\sharp}\), the following identification between quadratic forms holds
\[
L^{\sharp}=\Gamma^{\perp}/\Gamma
\]
 where the quadratic form on the right hand side is the quadratic form induced by \(M^{\sharp} \oplus N^{\sharp}\).
 We recall that two primitive embeddings \(M \hookrightarrow L\) are equivalent under the action of \(\OO(L)\) if and only if the corresponding groups \(H\) and \(H'\) are conjugate under the action of \(\OO(M)\) and \(\OO(N(-1))=\OO(N)\) in a way that commutes with the gluing isometries.
 
In an equivalent way, assuming that \(L\) is unique in its genus, we can give a primitive embedding \(M \hookrightarrow L\) giving a so called \textit{embedding subgroup} \(K \subset L^{\sharp}\), and an isometry \( \xi:K \to K' \subset M^{\sharp}\). We denote by \(\Sigma\) the graph of the isometry \(\xi\) in \(L^{\sharp} \oplus M(-1)^{\sharp}\). In this notation, the following identification between finite quadratic forms holds (similarly to what it is done for the gluing subgroup, the quadratic form on the right hand side is the one induced by \(L^{\sharp} \oplus M(-1)^{\sharp}\)):
\[
N^{\sharp}=\Sigma^{\perp}/\Sigma.
\]
Moreover if \(H\) is the gluing subgroup and if \(K\) is the embedding subgroup, the orders of these subgroups are called \textit{gluing index} and \textit{embedding index} respectively. These equalities hold:
\begin{equation}\begin{split}\label{formula per embedding}
h^{2} \cdot |\det(L)|=|\det(M) \cdot \det(N)|, \\ k^{2} \cdot |\det(N)|=|\det(L) \cdot \det(M)|.
\end{split}\end{equation}
In this paper, given a primitive embedding \(M \hookrightarrow L\) and given a vector \(v \in M\), we will need to compute the divisibility \((v,L)\) of the vector \(v\) in the whole lattice \(L\). 
\begin{remark}\label{relazione sulle divisibilità}
If \(M \hookrightarrow L\) is a primitive embedding and \(v\in M\), then it holds \((v,L) \ | \ (v,M)\).
\end{remark}

We recall here two results that will be useful in \autoref{Applications}.

\begin{proposition}\cite[Lemma 2.1]{gro_on_ven2021symplectic}\label{Lemma 2.1 GOV} If a primitive embedding of two lattices \(M \hookrightarrow L\) is defined by the gluing subgroup \(H \subset M^{\sharp}\), then it holds
\[
(v,L)= \mbox{\emph{max}} \{\ d \in \mathbb{N} \ | \ (v/d) \in H^{\perp}\}
\]
for every \(v \in M\). 
\end{proposition}

\begin{proposition}\cite[Corollary 2.2]{gro_on_ven2021symplectic}\label{cor 2.2 GOV} If \(M \hookrightarrow L\) is a primitive embedding and if \( |\det(M^{\perp})| = |\det(L)\cdot \det(M)|\), then \((v,L) = 1\) for every \(v \in M\).
\end{proposition}

\begin{proposition}\label{embedding in a unim lattice}
If \(M \hookrightarrow L\) is a primitive embedding and \(L\) is a unimodular lattice and if we denote by \(N=M^{\perp}\) the orthogonal complement of \(M\) in \(L\), then \(M^{\sharp} = N^{\sharp}\) as groups.
\end{proposition}

\proof
By \eqref{formula per embedding} it holds that \(
|H|^{2} \cdot |L^{\sharp}|=|M^{\sharp}| \cdot |N^{\sharp}|.
\) 
Since \(L\) is unimodular, then \begin{equation}\label{ref}
|H|^{2}=|M^{\sharp}| \cdot |N^{\sharp}|.
\end{equation}
There exists a  gluing isometry \(\gamma \colon H \to H'\) where \(H \subseteq M^{\sharp}\) and \(H' \subseteq N(-1)^{\sharp}\) hence \(|H|=|H'|\) and, in particular, \(|H| \leq |M^{\sharp}|\) and \(|H'| = |H|\leq |N^{\sharp}|\). By \eqref{ref} it holds that \(|H|=|M^{\sharp}|=|N^{\sharp}|\) and there exists an isometry of finite quadratic forms \(\gamma \colon M^{\sharp} \to N(-1)^{\sharp}\) hence, in particular, \(M^{\sharp}=N^{\sharp}\) as groups.
\endproof

Whenever we have a lattice \(L\) we can consider a subgroup \(G \subset \OO(L)\) of isometries of \(L\). We denote by \(L^{G}\) the \textit{invariant lattice} and its orthogonal complement \(L_{G}=(L^{G})^{\perp} \subset L\) is called \textit{coinvariant lattice}.

\subsubsection{The lattice structure of irreducible holomorphic symplectic manifolds}\label{The lattice structure of O'Grady's sixfolds}
If \(X\) is an irreducible holomorphic symplectic manifold, then the second integral cohomology is a torsion-free \(\mathbb{Z}\)-module of finite rank. Moreover there exists an integral symmetric bilinear form on \(\HH^{2}(X,\mathbb{Z})\) that endow the latter with a lattice structure. More precisely it holds the following theorem due to Beauville \cite{beauville1983varietes} and Fujiki  \cite{Fujiki_87_on_the_de_rham}.
\begin{theorem}
	Let \(X\) be a \(2n\)-dimensional irreducible holomorphic symplectic manifold. There exists a unique bilinear integral symmetric form \((-, -)_{X}\) defined on \(\HH^{2}(X,\mathbb{Z})\), the Beauville--Bogomolov--Fujiki form, and a unique positive constant \(c_{X}\), the Fujiki constant, such that for any \(\alpha \in \HH^{2}(X,\mathbb{Z})\) the following equality holds:
	\[
	\int_{X} \alpha^{2n}=c_X
(\alpha, \alpha)_{X}^{n},
	\]
	and for \(0 \neq \omega \in \HH^{0}(X,\Omega_X^{2})\)
	\[
	(\omega+\overline{\omega},\omega + \overline{\omega})>0.
	\]
\end{theorem}
\begin{remark}
	The Beauville--Bogomolow--Fujiki form \((-,-)_{X}\) and the Fujiki constant \(c_{X}\) are invariant up to deformation. 
\end{remark}

If \(X\) is an irreducible holomorphic symplectic manifold of \(\OG_6\) type, we denote by \(\bL\) the isometry class of \(\HH^{2}(X,\mathbb{Z})\), which depends only on the deformation type of \(X\). Rapagnetta \cite[Corollary 3.5.13]{rapagnetta2007topological} proves that
\[
\bL= \bU^{\oplus3} \oplus [-2]^{\oplus 2}.
\]

The second integral cohomology lattice \(\HH^{2}(X,\mathbb{Z})\) of an irreducible holomorphic symplectic manifold \(X\) is endowed with a pure weight-two Hodge structure. According to \autoref{Lattices_SUB}, if \((X, \eta)\) is a marked pair where \(X\) is an irreducible holomorphic symplectic manifold of \(\OG_6\) type then 
\[
\bL^{1,1}=\bL_{\mathbb{C}}^{1,1} \cap \bL,
\]
The integral lattice \(\bL^{1,1}\) is the Néron–Severi lattice of the marked pair \((X, \eta)\) of \(\OG_6\) type.

\subsection{O'Grady's sixfolds}\label{O'Grady's sixfold}
So far three deformation families of irreducible holomorphic symplectic manifolds in dimension six are known: manifolds of K3\(^{[3]}\) type, manifolds of Kum\(_n(\A)\) type and manifolds of \(\OG_6\) type. The latter class was discovered by O'Grady \cite{o2000new} as a resolution of singularities of a moduli space of sheaves on an abelian surface \(\A\).

\subsubsection{The Mukai lattice}\label{The Mukai lattice}
If \(\A\) is an abelian surface, we denote by \(\widetilde{\HH}(\A,\mathbb{Z})\) the even integral cohomology of \(\A\) i.e.
\[
\widetilde{\HH}(\A,\mathbb{Z})=\HH^{2*}(\A,\mathbb{Z})=\HH^{0}(\A, \mathbb{Z}) \oplus \HH^{2}(\A, \mathbb{Z}) \oplus \HH^{4}(\A, \mathbb{Z}).
\]
The \(\mathbb{Z}\)-module \(\widetilde{\HH}(\A, \mathbb{Z})\) has a lattice structure due to a pairing defined on it, the Mukai's pairing given by:
\[
(r_1,l_1,s_1)(r_2,l_2,s_2)=l_1l_2-r_1s_2-r_2s_1.
\]
where \(r_i \in \HH^{0}\), \(l_i \in \HH^{2}\) and \(s_i \in \HH^4\). This lattice is referred to as the \textit{Mukai lattice} of \(\A\) and it is isometric to \(\bU^{\oplus4}\). An element \(v \in \widetilde{\HH}(\A,\mathbb{Z})\) will be written as \((v_0,v_1,v_2)\), and if \(v_0 \geq 0\) and \(v_1 \in \NS(\A)\) then \(v\) is called \textit{Mukai vector}.

Moreover \(\widetilde{\HH}(\A, \mathbb{Z})\) has a pure weight-two Hodge structure such that the \((2,0)\)-part and the \((0,2)\)-part of \(\widetilde{\HH}(\A, \mathbb{C})\) are \(\HH^{2,0}(\A)\) and \(\HH^{0,2}(\A)\) respectively and the \((1,1)\)-part concerns of the following contributes:
\[
\widetilde{\HH}^{1,1}(\A)=\HH^{0}(\A, \mathbb{C}) \oplus \HH^{1,1}(\A) \oplus \HH^{4}(\A,\mathbb{C}).
\]
If \(v \in \widetilde{\HH}(\A,\mathbb{Z})\) is a Mukai vector then the sublattice with respect to the Mukai pairing

\begin{equation}\label{vperp}
v^{\perp}=\{ \alpha \in \widetilde{\HH}(\A,\mathbb{Z}) \ | \ (\alpha, v)=0\} \subseteq \widetilde{\HH}(\A,\mathbb{Z})
\end{equation}
inherits a pure weight-two Hodge structure from the one on \(\widetilde{\HH}(\A,\mathbb{Z})\).
More precisely it holds that
\begin{equation}\begin{split}\label{lattice structure of vperp}
(v^{\perp})^{0,2}=(v^{\perp} \otimes \mathbb{C}) \cap \widetilde{\HH}^{0,2}(\A),\\
(v^{\perp})^{2,0}=(v^{\perp} \otimes \mathbb{C}) \cap \widetilde{\HH}^{2,0}(\A),\\(v^{\perp})^{1,1}=(v^{\perp} \otimes \mathbb{C}) \cap \widetilde{\HH}^{1,1}(\A),
\end{split}\end{equation}

If \(\mathscr{F}\) is a coherent sheaf on \(\A\), its \textit{Mukai vector} is defined as follows:
\[
v(\mathscr F) = \Ch(\mathscr F) \sqrt{td(A)}=(\rank(\mathscr F), c_1(\mathscr F), ch_2(\mathscr F)) \in \widetilde{\HH}(\A,\mathbb{Z}).
\]
By construction for any coherent sheaf \(\mathscr F\) its Mukai vector is of \((1,1)\)-type and it satisfies one of the following relations

\begin{itemize}
    \item \(r>0\),
    \item \(r=0\) and \(l \neq 0\) with \(l\) effective,
    \item \(r=l=0\), and \(s>0\).
\end{itemize}
By \cite[Definition 2.27]{MW14} we have the following definition.

\begin{definition}\label{positive Mukai vector}
A vector \(v \in \widetilde{\HH}(\A, \mathbb{Z})\), \(v \neq 0\) satisfying \(v^2 \geq 2\) and the conditions above is called a \textit{positive Mukai vector}.
\end{definition}

\subsubsection{Moduli spaces of sheaves of \(\OG_6\) type}\label{O'Grady's sixfold starting from OLS-triple}
Let \(\theta\) be a \(v\)-generic polarization and \(v\) a  Mukai vector on \(\A\). We write \(M_{v}(\A, \theta)\) (respectively \(M_{v}^{s}(\A, \theta)\)) for the moduli space of \(\theta\)-semistable (resp \(\theta\)-stable) sheaves on the abelian surface \(\A\), with Mukai vector \(v\). We consider a Mukai vector \(v=mw\) where \(m \in \mathbb{N}\) and \(w\) is a primitive Mukai vector on \(\A\). It is well known that, if \(M_{v}^{s}(\A,\theta) \neq \emptyset\), then \(M_{v}^{s}(\A,\theta)\) is smooth of dimension \(v^{2}+2\) and carries a symplectic form (see Mukai \cite{mukai1984symplectic} for more details). Since we are taking into consideration a moduli space on an abelian surface, a further construction is necessary: choose \(\mathscr F_{0} \in M_{v}(\A, \theta)\), and define the following map \cite{yoshioka1999some}: 
\begin{equation}\label{map for OG6}
a_{v}: M_{v}(\A, \theta) \longrightarrow \A\times \A^{\vee}
\end{equation}
\[
a_{v}(\mathscr F):=(det(p_{\A^{\vee}!}((\mathscr F -\mathscr F_{0}) \otimes (\mathscr P - \mathscr O_{\A \times \A^{\vee}})), det(\mathscr F) \otimes \det(\mathscr F_{0})^{-1}),
\]
where \(p_{\A^{\vee}}\colon \A \times \A^{\vee} \longrightarrow {\A}^{\vee}\) is the projection and \(\mathscr P\) is the Poincaré bundle on \(A\times A^{\vee}\).

We define
\[\K_{v}(\A, \theta)= a_{v}^{-1}(0_{\A}, \mathscr O_{\A}),\]
where \(0_{\A}\) is the zero of \(\A\). 
We recall the following crucial result in the case \(v\) is a primitive Mukai vector:

\begin{theorem}\cite[Theorem 0.2]{yoshioka2001moduli}
	Let \(\A\) be an abelian surface and let \(v\) be a primitive Mukai vector, let \(\theta\) be a \(v\)-generic polarization. Then \(M_{v}(\A, \theta)=M_{v}^{s}(\A, \theta)\). If \(v^{2} \geq 6\) then \(\K_{v}(\A, \theta)\) is an irreducible holomorphic symplectic  manifold of dimension \(2n=v^{2}-2\), which is deformation equivalent to \(Kum_{n}(\A)\), the generalized Kummer variety of \(\A\), and there is a Hodge isometry between \(v^{\perp}\) (see \eqref{vperp}) and \(\HH^{2}(\K_v(\A,\theta), \mathbb{Z})\).
\end{theorem}

If the Mukai vector \(v\) is not primitive, then \(M_v(\A,\theta)\) can be singular. O'Grady started from this consideration to find a new deformation class of irreducible holomorphic symplectic manifolds.
It is natural to ask if there is a symplectic resolution of singularities \(\pi_v \colon \widetilde{M_v}(\A,\theta) \to  M_v(\A,\theta)\), such that on \(\widetilde{M_v}(\A,\theta)\) there is a symplectic form extending the one on \(M_v^s(\A,\theta)\). The first result appearing in literature is the one of O'Grady. For more details see \cite{o2000new}.
\begin{theorem}\cite[Theorem 1.4]{o2000new} Let \(\A\) be an abelian surface, \(v=(2,0,-2)\) and \(\theta\) a \(v\)-generic polarization. Then \(\K_6=\K_v(\A, \theta)\) admits a symplectic resolution \(\pi:\widetilde{\K}_6 \to \K_6\) and \(\widetilde{\K}_6\) is an irreducible symplectic variety of dimension \(6\) and second Betti number \(8\). Manifolds deformation equivalent to \(\widetilde{\K}_6\) are called manifolds of \(\OG_6\) type.
\end{theorem}

Moreover what is done by O'Grady for a specific Mukai vector was generalized by Perego--Rapagnetta for a more general class of surfaces, Mukai vectors and polarizations. More precisely, Perego--Rapagnetta introduce the OLS-triple, after the work of O'Grady\cite{o1997desingularized}\cite{o2000new}, and Lehn--Sorger \cite{Lehn_Sorger}(see \cite[Definition 1.5]{perego2013deformation} for the definition of OLS-triple), to find the more general setting of Mukai vectors and polarizations that admit an analogue result of the one of O'Grady. If \((\A, v, \theta)\) is an OLS-triple, then \(M_v(\A, \theta)\) admits a symplectic resolution \(\widetilde{M}_v(\A,\theta)\) obtained as the blow-up of \(M_v(\A,\theta)\) along the singular locus \(\Sigma_v=M_v(\A, \theta) \setminus M_v^s(\A,\theta)\) with reduced structure. Moreover 
\[
\widetilde{\K}_v(\A, \theta)=\widetilde{\K}_v(\A,\theta)=\pi_v^{-1}(\K_v(\A,\theta)),
\]
and we still write \(\pi_v:\widetilde{\K}_v(\A,\theta) \to \K_v(\A,\theta)\) for the symplectic resolution.
They give the following result.
\begin{theorem}\cite[Theorem 1.6]{perego2013deformation}
Let \((\A, v, H)\) be an OLS-triple where \(\A\) is an abelian surface. The moduli space \(\widetilde{\K}_v(\A, \theta)\) is an irreducible holomorphic symplectic manifold which is deformation equivalent to \(\widetilde{\K}_6\).
\end{theorem}
Moreover Perego--Rapagnetta gives a result about the weight-two Hodge structure of the second integral cohomology of \(\widetilde{\K}_v(\A,\theta)\).

\begin{theorem}\cite[Theorem 1.7]{perego2013deformation} Let \(\A\) be an abelian surface and let \((\A,v,\theta)\) be an OLS-triple. The pullback \(\pi^*:\HH^{2}(\K_v(\A,\theta), \mathbb{Z}) \to \HH^{2}(\widetilde{\K}_v(\A,\theta), \mathbb{Z})\) is injective, and the
restrictions to \(\HH^2(\K_v(\A,\theta),\mathbb{Z})\) of the pure weight-two Hodge structure and of the Beauville--Bogomolov--Fujiki form on \(\HH^2(\widetilde{\K}_v(\A,\theta),\mathbb{Z})\) give a pure weight-two Hodge structure on \(\HH^{2}(\K_v(\A,\theta), \mathbb{Z})\) and a lattice structure on \(\HH^2(\K_v(\A,\theta),\mathbb{Z})\). Moreover, there is an isometry of weight-two Hodge structures
\[
\nu_v \colon \HH^{2}(\K_v(\A,\theta), \mathbb{Z}) \xrightarrow{\sim} v^{\perp} \subset \widetilde{\HH}(\A,\mathbb{Z})
\]
\end{theorem}

Moreover Perego--Rapagnetta \cite{Per_Rap_factoriality} computed the lattice and Hodge structure of \(\HH^{2}(\widetilde{\K}_{v}(\A, \theta), \mathbb{Z})\) in terms of the Hodge structure of \(v^{\perp}\) as a sublattice of the Mukai lattice \(\widetilde{\HH}(\A,\mathbb{Z})\) introduced in \autoref{The Mukai lattice}. They consider the \(\mathbb{Z}\)--module \(v^{\perp} \oplus_{\perp} \mathbb{Z} \cdot\sigma\) with the symmetric bilinear form on \(v^{\perp}\) induced by the Mukai pairing and defining  \((\sigma,\sigma )=-2\).
The \(\mathbb{Z}\)--module \(v^{\perp} \oplus_{\perp} \mathbb{Z} \cdot\sigma\) is a lattice and carries a pure weight-\(2\) Hodge structure in the following way:

\[
 (v^{\perp} \oplus_{\perp} \mathbb{Z} \cdot\sigma)^{2,0}=(v^{\perp})^{2,0} = \widetilde{\HH}^{2,0}(\A)=\HH^{2,0}(\A)
\]
\[
(v^{\perp} \oplus_{\perp} \mathbb{Z} \cdot\sigma)^{0,2}=(v^{\perp})^{0,2} = \widetilde{\HH}^{0,2}(\A)=\HH^{0,2}(\A)
\]
\[
(v^{\perp} \oplus_{\perp} \mathbb{Z} \cdot\sigma)^{1,1}=(v^{\perp})^{1,1} \oplus \mathbb{C} \cdot \sigma=\HH^{0}(\A) \oplus \HH^{1,1}(\A) \oplus \HH^4(\A) \oplus \mathbb{Z} \cdot \sigma 
\]

Note that \(\sigma^{\perp} \subset (v^{\perp} \oplus_{\perp} \mathbb{Z} \cdot\sigma)\cong v^{\perp} \subset \widetilde{\HH}(\A,\mathbb{Z})\). Note that the Hodge structure on \(v^{\perp}\) is the one recalled in section \autoref{The Mukai lattice}. Moreover note that the divisibility of \(\sigma\) in the lattice \((v^{\perp} \oplus_{\perp} \mathbb{Z} \cdot\sigma)\) is \(2\). 
\begin{theorem}\cite[Theorem 3.4]{Per_Rap_factoriality}\label{Hodge structure of Ktilde}
	There is a Hodge isometry of pure wight-two Hodge structures 
	\[
	v^{\perp} \oplus_{\perp} \mathbb{Z} \cdot \sigma \cong \HH^{2}( \widetilde{\K}_{v}(\A,\theta), \mathbb{Z})
	\]
   where the lattice structure on the right-hand side is given by the Beauville--Bogomolov--Fujiki quadratic form.
\end{theorem}

\section{Induced groups of automorphisms}\label{Induced automorphisms groups}
An example of irreducible holomorphic symplectic manifolds that arise from symplectic surfaces are the Hilbert schemes of \(n\) points on K\(3\) surfaces, constructed by Beauville in \cite{beauville1983varietes}. This kind of construction allows to produce several examples of automorphisms of irreducible symplectic manifolds of K3\(^{[n]}\) type, simply by taking a K3 surface with a non-trivial automorphism group and considering the induced action on its Hilbert scheme. These kinds of automorphisms are called \textit{natural} and were studied by Beauville \cite{beauville1983some}, Boissière \cite{boissiere2012automorphismes} and many others.
A natural question is to ask when a birational transformation of a manifold of \(\OG_6\) type, that is a moduli space of sheaves on an abelian surface, is induced by an automorphism of the abelian surface. In \autoref{proof of Thm 1.1} we prove a result to determine when a manifold \(X\) of \(\OG_6\) type is birational to a moduli space, and in \autoref{section_4.1} we prove a numerical criterion to determine when a birational transformation is induced.
%
%
%
%
%
%
\subsection{Proof of \autoref{X n.m.s allora birazionale alla fibra della mappa di albanese}}\label{proof of Thm 1.1}
This section is devoted to determine when a manifold \(X\) of \(\OG_6\) type is birational to \(\widetilde{\K}_v(\A,\theta)\) where \((\A,v,\theta)\) is an OLS-triple and whose construction is recalled in \autoref{O'Grady's sixfold starting from OLS-triple}.
We state a necessary and sufficient criterion entirely in terms of the lattice structure of the second integral cohomology of \(X\). In the following \(\bLambda_{8} = \bU^{\oplus 4}\) and \(\bLambda_{10} = \bU^{\oplus 5}\).

\begin{definition}\label{Hodge embedding}
Let \(L\) be a lattice endowed with a pure weight-two Hodge structure and consider the following primitive embedding in a lattice \(\Lambda\)
\[
i : L \hookrightarrow \Lambda.
\]
We call the embedding \(i\) a \textit{Hodge embedding} if \(\Lambda\) is endowed with a pure weight-two Hodge structure inherited by \(L\), defined as follows:

 \[
 \Lambda^{2,0}=L^{2,0},
\]
\[
\Lambda^{0,2}=L^{0,2},
\]
\[
\Lambda^{1,1}=L^{1,1} \oplus L^{\perp_{\Lambda}}.
 \]

\end{definition}

If \((X,\eta)\) is a marked pair of \(\OG_6\) type and if there exists a class \(\sigma \in \bL^{1,1}\) such that \(\sigma^{2}=-2\) and \((\sigma, \bL)=2\), then there is a unique \cite[Proposition 1.5.1]{Nik79} (up to isometry) primitive embedding \( \sigma \hookrightarrow \bL\) and \(\sigma^{\perp}\) inherits a pure weight-two Hodge structure in the following way:
\begin{equation}\begin{split}
(\sigma^{\perp})^{2,0}=\bL^{2,0}\\
(\sigma^{\perp})^{0,2}=\bL^{0,2}\\ (\sigma^{\perp})^{1,1}=\bL^{1,1}\cap \sigma^{\perp}\\
\end{split}\end{equation}
where the isometry class of \(\sigma^{\perp}\) is \(\bU^{\oplus3} \oplus [-2]\).
\begin{definition}\label{nummodspace}
	Let \((X, \eta)\) be a projective marked pair of \(\OG_6\) type, where \(\eta:\HH^{2}(X,\mathbb{Z}) \to \bL\) is a fixed marking. We call \(X\) a \textit{numerical moduli space} if there exists a class \(\sigma \in \bL^{1,1}\) such that \(\sigma^{2}=-2\) and \((\sigma, \bL)=2\) and through the Hodge embedding \( \sigma^{\perp} \hookrightarrow \bLambda_8\) the lattice \(\bLambda_{8}^{1,1}\) contains a copy of \(\bU\) as a direct summand.
\end{definition}

\begin{proposition}\label{relazione segnature}
Let \((X,\eta)\) be a projective marked pair of \(\OG_6\) type. If \(X\) is a numerical moduli space then 
\[
\sign(\bLambda_8^{1,1})=\sign(\bL^{1,1})+(1,-1)
\]
\end{proposition}
\proof
By assumption there exists a negative class \(\sigma \in \bL^{1,1}\). We denote by \([\sigma]\) the lattice of rank \(1\) and signature \((0,1)\) generated by \(\sigma\). By construction it holds that \[
\sign(\bL^{1,1})=\sign((\sigma^{\perp})^{1,1}) + \sign([\sigma]).
\] 
The orthogonal complement of \(\sigma^{\perp} \cong \bU^{\oplus3}\oplus[-2]\) in \(\bLambda_8\) is a rank \(1\) lattice of signature \((1,0)\) that we denote by \([w]\). The class \(w \in \bLambda_8^{1,1}\) is of \((1,1)\) type by definition of Hodge embedding, hence it holds the following equality
\[
\sign(\bLambda_8^{1,1})=\sign((\sigma^{\perp})^{1,1}) + \sign([w]).
\]
From the two previous relations we have 
\[
\sign(\bLambda_8^{1,1})=\sign(\bL^{1,1}) - \sign([\sigma]) + \sign([w])= \sign(\bL^{1,1}) + (1,-1).
\]\endproof

 In the following we denote by \(\widetilde{\K}_v(\A,\theta)\) an irreducible holomorphic symplectic manifold of \(\OG_6\) type obtained starting from an OLS-triple \((\A,v,\theta)\), as recalled in section \autoref{O'Grady's sixfold starting from OLS-triple}.

\proof (of \autoref{X n.m.s allora birazionale alla fibra della mappa di albanese})
If \(\Phi \colon X \dashrightarrow \widetilde{\K}_v(\A, \theta)\) is a birational morphism, then the induced isometry \(\Phi^{*}:\HH^{2}(\widetilde{\K}_v(\A, \theta), \mathbb{Z}) \to \HH^{2}(X,\mathbb{Z})\) is an isometry of Hodge structures. The manifold \(\widetilde{\K}_v(\A, \theta)\) is obtained as a resolution of the singular moduli space \({\K}_v(\A, \theta)\). The class of the exceptional divisor \(E \in \HH^2(\widetilde{\K}_v, \mathbb{Z})\) is of \((1,1)\) type, of square \(-2\) and divisibility \(2\), and \(E^{\perp} \cong \HH^{2}(\K_v, \mathbb{Z})\) \cite[Theorem 3.4(2)]{Per_Rap_factoriality}. Moreover there exists a Hodge embedding (see \autoref{Hodge embedding}) \(\HH^{2}(\K_v, \mathbb{Z}) \hookrightarrow \widetilde{\HH}(\A,\mathbb{Z}) \cong \bLambda_8\). By construction, the induced weight-two Hodge structure on \(\widetilde{\HH}(\A,\mathbb{Z})\) is the one defined in \autoref{The Mukai lattice} on the Mukai lattice. The vectors \((1,0,0)\) and \((0,0,1)\) generates a copy of \(\bU\) in the lattice \(\widetilde{\HH}^{1,1}(\A, \mathbb{Z})=\bLambda_8^{1,1}\), hence \(X\) is a numerical moduli space.

For the other direction, let \( \HH^2(X,\mathbb{Z}) \to \bL\) be a fixed isometry, by assumption there exists a class \(\sigma \in \bL^{1,1}\) such that \(\sigma^2=-2\), \((\sigma, \bL)=2\), and a Hodge embedding \(\sigma^{\perp}\hookrightarrow \bLambda_8\)
such that \(\bLambda_8^{1,1}\) contains \(\bU\) as a direct summand. Denote by \(w\) the orthogonal complement of \(\sigma^{\perp}\) in \(\bLambda_8\), then it holds that \(w^2=2\) and \(w^{\perp}(\subset \bLambda_8) \cong \bU^{\oplus3} \oplus [-2]\). The lattice \(w^{\perp}\) is a sublattice of \(\bLambda_8\) hence it inherits a weight-two Hodge structure. Note that \(w^{\perp}(\subset \bLambda_8) \cong \sigma^{\perp}(\subset \bL)\). Moreover the signature of \(\sigma^{\perp}\) is equal to \((3,4)\) and since \(X\) is a numerical moduli space, \(X\) is projective, hence the positive part of the signature of \(\NS(X) = \bL^{1,1}\) is equal to \(1\). The class \(\sigma\) has negative square hence the positive part of the signature of \((\sigma^{\perp})^{1,1}\) is equal to the positive part of the signature of \(\bL^{1,1}\). Thus we get \((\bLambda_8)^{1,1}=(\sigma^{\perp})^{1,1} \oplus [w]\) hence the positive part of the signature of \((\bLambda_8)^{1,1}\) is equal to \(2\). Moreover the Hodge structure induced on \(\bLambda_8\) is the one inherited from \(\sigma^{\perp} \subset \bL\) hence \(\bL^{2,0}=(\bLambda_8)^{2,0}\), \(\bL^{0,2}=(\bLambda_8)^{0,2}\), and consequently the signature of \((\bLambda_{8})^{1,1}\) is equal to \((2,4)\). We have \((\bLambda_8)^{1,1}=\bU \oplus \T\) where \(\T\) is an even lattice of signature \((1,3)\). By \cite[Theorem 2.4]{mongardi2015prime_abelian} there exists an abelian surface \(\A\) such that \(\T=\NS(\A)\) and we call \(\theta\) the generator of the positive part of \(\T\). We define the Mukai vector \(v=2w\) and we can choose a \(v\)-generic polarization \(\theta\) on \(\A\). Now by \cite[Lemma 2.28]{MW15} since \(w^2=2\) then \(w\) or \(-w\) is a positive Mukai vector ( see \autoref{positive Mukai vector}) hence we may assume the \((\A,v,\theta)\) is an OLS-triple. Then we consider the map \(a_v \colon M_v(\A, \theta) \to \A \times \A^{\vee}\) and we define \(\K_v(\A, \theta)=a_v^{-1} (0, \mathcal{O}_{\A})\). By \cite[Theorem 1.7]{perego2013deformation} there exists a Hodge isometry 
\begin{equation}
\HH^{2}(\K_{v}(\A,\theta), \mathbb{Z}) \xrightarrow{\sim} v^{\perp} \subset \widetilde{\HH}(\A, \mathbb{Z}) \cong \bLambda_{8},
\end{equation}
where the orthogonal complement of \(v \in \widetilde{\HH}(\A, \mathbb{Z})\) is computed with respect to the Mukai pairing.
 The singular moduli space \(\K_{v}(\A,\theta)\) admits a symplectic crepant resolution \(\widetilde{\K}_{v}(\A,\theta)\) and the exceptional divisor \(E\) is such that \(E^{2}=-2\) \cite[Corollary~3.5.13]{rapagnetta2007topological}. It holds the following Hodge isometry \(\HH^{2}(\widetilde{\K}_{v}(\A,\theta), \mathbb{Z}) \cong \HH^{2}(\K_{v}(\A,\theta), \mathbb{Z}) \oplus \mathbb{Z} \cdot E\).
 Finally we obtain the Hodge isometry
 \[
 \sigma^{\perp} \oplus \sigma \xrightarrow{\sim} w^{\perp} \oplus \mathbb{Z} \cdot E
 \]
 which means that there exists a Hodge isometry between \(\bL=\sigma^{\perp} \oplus \sigma\) and \(\HH^{2}(\widetilde{\K}_v(\A, \theta), \mathbb{Z})\)
 hence between \(\HH^{2}(X,\mathbb{Z})\) and \(\HH^{2}(\widetilde{\K}_v(\A, \theta), \mathbb{Z})\). By \cite[Theorem 5.2(2)]{MROG6mon} \(X\) is birational to \(\widetilde{\K}_{v}(\A,\theta)\).
\endproof

\subsection{Induced birational transformations and proof of \autoref{num ind impliesd ind}}\label{section_4.1}
In this section we give a lattice-theoretic criterion to determine when a birational transformation of a manifold \(X\) of \(\OG_6\) type is induced by an automorphism of the abelian surface \(\A\), in the case in which \(X\) is birational to the moduli space \(\widetilde{\K}_v(\A,\theta)\).

Consider an automorphism of the abelian surface \(\A\). It induces an isometry of \(\widetilde{\HH}(\A, \mathbb{Z})\). Moreover if the automorphism fixes the Mukai vector \(v\) then we obtain an isometry of \( \HH^{2}( \widetilde{\K}_{v}(\A,\theta), \mathbb{Z}) \) asking that \(\sigma\) is fixed.

In the following we give all the statements assuming that \(\G \subset \Bir(X)\) is a group of birational transformations of \(X\). If the statements holds true for \(\G \subset \Aut(X)\) then \(\G\) will be called an \textit{induced group of automorphisms} or a \textit{numerically induced group of automorphisms}.

\begin{definition}\label{numerically induced}
	Let \(X\) be a smooth projective irreducible holomorphic symplectic manifold of \(\OG_6\) type and let \(\eta:\HH^{2}(X,\mathbb{Z}) \to \bL\) be a marking. Let \(\G \subset \Bir(X)\) be a finite group of birational transformations. Assume that there exists a class \(\sigma\) of \((1,1)\) type on \(X\) such that \(\sigma^2=-2\) and \((\sigma, \bL)=2\), and consider the primitive \textit{Hodge embedding} \(i:\sigma^{\perp} \hookrightarrow \bLambda_8\). 
	The group \(\G \subset \Bir(X)\) is called a \textit{numerically induced group of birational transformations} if the following hold:
	\begin{itemize}
		\item[(1)] the class $\sigma \in \NS(X)$ is \(\G\)-invariant;
		\item [(2)] the induced action of \(\G\) on \(\bLambda_{8}\) is such that the \((1,1)\) part of the invariant lattice \((\bLambda_8)^{\G}\) contains $\bU$ as a direct summand;
		\item[(3)] 	for all $g \in \G$, $\det(g^*)=1$.
	\end{itemize}
\end{definition}

\begin{proposition}\label{aut di A induced aut mod space}
If \(\varphi \in \Aut(\A)\) is an automorphism of the abelian surface \(\A\), \(v \in \widetilde{\HH}(\A,\mathbb{Z})\) is a \(\varphi\)-invariant Mukai vector on \(\A\), and \(\theta\) is a \(\varphi\)-invariant polarization on \(\A\), then \(\varphi\) induces an automorphism on the fibre \(\widetilde{\K}_{v}(\A, \theta)\). Moreover the automorphism is numerically induced. 
\end{proposition}	
\proof
To prove that the automorphism \(\varphi\) of \(\A\) induces an automorphism of the moduli space \(M_{v}(\A, \theta)\) we need to check that the pullback along \(\varphi\) induces an automorphism of the moduli functor. From the definition of stability, if a sheaf \(\mathscr{F}\) is \(\theta\)-stable then \(\varphi^{*}\mathscr{F}\) is \(\theta\)-stable, see \cite[Proposition 2.32]{MW15}. Moreover if \(\varphi \in \Aut(\A)\) is an automorphism then by definition it preserves the origin hence the induced automorphism respects the fibre \(\K_{v}(\A,\theta)\) over \((0, \mathcal{O}_{\A})\) of the map \(a_v:M_{v}(\A, \theta) \to \A \times \A^{\vee}\) and we call \(\hat{\varphi}\) the induced action on it. 
Furthermore, the singular locus \(\Sigma\) is certainly \(\hat{\varphi}\)-invariant, and we have a well-defined induced
action on the normal bundle \(\mathcal{N}=\mathcal{N}_{\Sigma|\K_v(\A, \theta)}\). In fact the fibre \(\mathcal{N}_{\mathcal{F}_1,\mathcal{F}_{2}}\) of \(\mathcal{N}\) over \(\mathcal{F}_1 \oplus \mathcal{F}_2 \) is isomorphic to \( \Ext^{1}(\mathcal{F}_1, \mathcal{F}_2) \oplus \Ext^1(\mathcal{F}_2, \mathcal{F}_1)\) and we have that the map \(\mathcal{N}_{\mathcal{F}_1, \mathcal{F}_2} \to \mathcal{N}_{\varphi^{*}\mathcal{F}_{1}, \varphi^{*}\mathcal{F}_{2}}\) is the pullback map. Hence we get the induced action \(\widetilde{\varphi}\) on the blow up of the fibre \(\widetilde{\K}_v(\A, \theta)\). We call \(\sigma\) the class of the exceptional divisor which is obviously fixed by the induced action. Moreover, due to the primitive embedding \( \sigma^{\perp} \subset \HH^{2}(\widetilde{\K}_v(\A),\mathbb{Z}) \hookrightarrow \widetilde{\HH}(\A,\mathbb{Z}) \cong \bLambda_8\), the induced action of \(\widetilde{\varphi}\) on \(\sigma^{\perp}\) induces an action on \( \widetilde{\HH}(\A,\mathbb{Z}) \cong \bLambda_8\) which is the same action induced by \(\varphi\) on \(\widetilde{\HH}(\A,\mathbb{Z})\).
By assumption the sublattice \(\HH^{0}(\A, \mathbb{Z}) \oplus \HH^{4}(\A, \mathbb{Z})\) is contained in the \((1,1)\) part of the lattice \(\bLambda_8\). The class of the surface i.e. the generator of \(\HH^{4}(\A, \mathbb{Z})\), and the class of the points i.e.\ the generator of  \(\HH^{0}(\A, \mathbb{Z})\) are preserved by \(\varphi\), hence a copy of \(\bU\) is contained in the \((1,1)\) part of the lattice \((\bLambda_8)^{\varphi}\). 
Finally, since \(\varphi\) is an automorphism of \(\A\), the induced isometry on \(\HH^2(\A,\mathbb{Z})\) is, in particular, a monodromy operator hence by \cite[Section 3]{MROG6mon} it belongs to \(\SO^{+}(\HH^{2}(\A,\mathbb{Z}))\) i.e. its determinant is equal to \(1\).
\endproof
	
\begin{remark}\label{rem}
In the assumptions of \autoref{aut di A induced aut mod space} we ask that the polarization \(\theta \in \NS(\A)\) is \(\varphi\)-invariant, then the automorphism of the abelian surface induces an automorphism of the desingularized moduli space \( \widetilde{\K}_{v}(\A, \theta) \). This condition is never verified for a symplectic automorphism of \(\A\) except in the case in which the automorphism of \(\A\) is trivial. If \(\theta\) is not \(\varphi\)-invariant, then we get at least a birational self-map of \(\widetilde{\K}_v(\A, \theta)\). 
\end{remark}

\begin{definition}\label{induced group}
	Let \(X\) be a manifold of \(\OG_6\) type which is a numerical moduli space, and let \(\G \subset \Bir(X)\) be a finite subgroup of birational transformations of \(X\). The group \(\G \subset \Bir(X)\) is called an \textit{induced group of birational transformations} if there exists a group \(\G \subset \Aut(\A)\), there exists a Mukai vector \(v \in \widetilde{\HH}(\A, \mathbb{Z})^{\G}\) and a \(v\)-generic polarization \(\theta\) and the action induced by \(\G\) on \(\widetilde{\K}_{v}(\A, \theta)\) coincides with the given action of \(\G\) on \(X\) (up to automorphisms of \(\Ker(\eta_{*})\)).
\end{definition}

\begin{corollary}\label{G induced implies G num ind}
	Let \(X\) be a manifold of \(\OG_6\) type which is a numerical moduli space. Let \(\G \subset \Bir(X)\) be a finite subgroup. If \(\G\) is an induced group of birational transformations then \(\G\) is a numerically induced group of birational transformations.
\end{corollary}
\proof 
For the proof we refer to \autoref{aut di A induced aut mod space}. 
\endproof

\proof (of \autoref{num ind impliesd ind})
First of all let us consider the case in which the group \(\G\) is symplectic. Then we have \(\T(X) \subseteq \bL^{\G}\). Since \(\G\) is numerically induced then the class \(\sigma \in \NS(X)\) such that \(\sigma^2= -2\) and \((\sigma,\bL)=2\) is fixed by \(\G\), and we have the primitive Hodge embedding
\begin{equation}
 \sigma^{\perp} \hookrightarrow \bLambda_8.
\end{equation}
We call \(w\) the generator of the orthogonal complement of \(\sigma^{\perp}\) in \(\bLambda_8\) and we note that by construction \(w\) is fixed by the induced action of \(\G\) on \(\bLambda_8\).
Since \(\G\) is numerically induced then \((\bLambda_8)^{\G} = \bU \oplus \T\) and \(\bU\) is in the \((1,1)\) part of the Hodge structure of \((\bLambda_8)^{\G}\). Due to the fact that \(\sigma\) is \(\G\)-invariant then \(\bL_{\G} = (\sigma^{\perp})_{\G}= (\bLambda_8)_{\G}\). We then have that \(\bL_{\G}\) embeds in the abelian lattice and its orthogonal is \(\T\), where the action of \(\G\) is trivial. We give to this abelian lattice the induced Hodge structure from \(\bLambda_8\) and we let \(\A\) the corresponding abelian surface \cite[Theorem 2]{shioda1978period}. We can take \( v=2w \in \widetilde{\HH}(\A,\mathbb{Z}) \cong \bLambda_8\) as Mukai vector on \(\A\). By construction \(v\) is fixed by the induced action of \(\G\) on \(\bLambda_8\) and by \autoref{X n.m.s allora birazionale alla fibra della mappa di albanese} \(X\) is birational to the moduli space \(\widetilde{\K}_v(\A,\theta)\). We have the two following isometries of Hodge structures \begin{equation}\begin{split}
\HH^{2}(X,\mathbb{Z}) \to \HH^{2}(\widetilde{\K}_v(\A,\theta), \mathbb{Z})\\
\HH^{2}(\widetilde{\K}_v(\A,\theta), \mathbb{Z}) \to v^{\perp} \oplus_{\perp} \mathbb{Z} \cdot \sigma 
\end{split}\end{equation}
where the second one holds true by \autoref{Hodge structure of Ktilde}. In \autoref{O'Grady's sixfold starting from OLS-triple} we have described the pure weight-two Hodge structure on \(v^{\perp} \oplus \mathbb{Z} \cdot \sigma\) and the relation with the Hodge structure on \(\HH^{2}(\A,\mathbb{Z})\). An element of \(\G\) induces an isometry of Hodge structures on \(\HH^{2}(X,\mathbb{Z})\) hence composing the previous isometries of Hodge structures we have an isometry of Hodge structure on \(v^{\perp}\) hence also on \(\HH^{2}(\A,\mathbb{Z})\). By construction the group \(\G\) is a group of Hodge isometries on \(\A\), moreover \(\G\) is numerically induced, hence the isometries are orientation preserving and of determinant \(1\) hence by \cite[Section 3]{MROG6mon} they are in the monodromy group of \(\A\). Therefore by \cite[Theorem 2.1]{mongardi2015prime_abelian} the group \(\G\subset \Aut(\A)\) is a group of automorphism of \(\A\) and since by construction the Mukai vector \(v\) is preserved by \(\G\), by \autoref{aut di A induced aut mod space} and \autoref{rem} we have an induced group of birational transformations on \(\widetilde{\K}_v(\A, \theta)\). 
The induced action on the second integral cohomology of \(\widetilde{\K}_v(\A,\theta)\) is the action we started with, up to isomorphisms of the kernel of the representation map 
\[
\eta_{*} \colon \Bir(X) \to \OO(\bL),
\]
which is isomorphic to \(\A[2] \times \A^{\vee}[2]\) \cite[Theorem 5.2]{MW14}.

Now we assume that every non-trivial element of \(\G\) has a nonsymplectic action. In this hypothesis we can assume that \(\bL^{\G} = \NS(X)\). As in the previous case we have \((\bLambda_8)^{\G} = \bU \oplus \T\) and we can consider the abelian surface \(\A\) associated to the Hodge structure induced on \(\bU^{\perp} \subset \bLambda_8\). Again by \autoref{X n.m.s allora birazionale alla fibra della mappa di albanese} \(X\) is birational to a moduli space of sheaves on \(\A\). The group \(\G\) is a group of Hodge isometries of \(\A\) preserving \(\T=\NS(\A)\) and we can conclude as in the symplectic case.

Finally, if \(\G_s\) is the symplectic part of the group \(\G\), then we obtain an abelian surface \(\A\) as in the first step with \(\G_s\subset \Aut(\A)\). We can extend the action of the group \(\G\) on \(\A\) by applying the second step to the quotient group \(\hat{\G} = \G/\G_s\).
\endproof
\begin{remark}
There exist automorphisms of manifolds of \(\OG_6\) type that act trivially on the second integral cohomology \cite[Theorem 5.2]{MW14} and this explains the factor  \(\A^{\vee}[2]\) in the theorem above. 
\end{remark}

\begin{corollary} \label{S_G(X) should be even in tha p=2 case}
	Let \((X, \eta)\) be a marked irreducible holomorphic symplectic manifold of \(\OG_6\) type and let \(\eta \colon \HH^{2}(X,\mathbb{Z}) \to \bL\) be a marking. If the group \(\G \subset \Bir(X)\) is an induced group of birational transformations and \(|\G|=2\), then \(\rank(\bL_{\G})\) is even.
\end{corollary}
\proof
If \(\G \subset \Bir(X)\) is induced then by \autoref{G induced implies G num ind} it is numerically induced hence by \autoref{numerically induced} every element in \(\G\) has an induced action in cohomology with determinant \(1\). If \(|\G|=2\) and \(\varphi\) is a generator of \(\G\) then \(\det(\varphi)=(-1)^{\rank(\bL_{\G})}\) and this implies that \(\rank(\bL_{\G})\) is even.
\endproof


\section{Automorphisms induced at the quotient} \label{Automorphisms induced at the quotient}

In this section we refer to the construction of the birational model of manifolds of \(\OG_6\) type described by Mongardi--Rapagnetta--Sacc\`{a} \cite{MRSHodge} (briefly MRS construction). In \autoref{The Mongardi--Rapagnetta--Sacca model} we recall the main steps of the MRS construction which provides a birational model of manifolds of \(\OG_6\) type as the quotient of a manifold of K3\(^{[3]}\) type by a birational symplectic involution. In \autoref{proof thm 1.3} we consider a manifold \(X\) of \(\OG_6\) type that is birational to a moduli space \(\widetilde{\K}_v(\A,\theta)\), and we prove a result to determine when a birational transformation of \(X\) lifts to a birational transformation of the manifold of K3\(^{[3]}\) type involved in the MRS construction. Finally in \autoref{A sufficient cond to have a regular morphism of K3} we prove \autoref{thm induced at the quotient} to determine when the birational transformation of \(X\) lifts to a regular automorphism of the manifold of K3\(^{[3]}\) type.

\subsection{The Mongardi--Rapagnetta--Sacc\`{a} model}\label{The Mongardi--Rapagnetta--Sacca model}
Mongardi--Rapagnetta-Sacc\`{a} show that for any abelian surface \(\A\), for an effective Mukai vector  (the Mukai vector of a coherent sheaf on \(\A\)) \(v=2v_0\) with \(v_0^2=2\) on \(\A\), and for a \(v\)-generic principal polarization \(\theta\) on \(\A\), the irreducible holomorphic symplectic manifold of dimension six \(\widetilde{\K}_{v}(\A, \theta)\) admits a rational double cover from a normal projective variety which is birational to an irreducible holomorphic symplectic manifold \(\underline{Y}_v(\A,\theta)\) of K3\(^{[3]}\) type.
More precisely, the singular locus \(\Sigma_v \subset \K_{v}(\A, \theta)\) has codimension 2. The inverse image \(\widetilde{\Sigma}_v\) of \(\Sigma_v\) in \(\widetilde{\K}_v(\A,\theta)\) is an irreducible divisor, which is divisible by two in the integral cohomology by a result of Rapagnetta \cite[Theorem 3.3.1]{rapagnetta2007topological}. The Picard group of an irreducible holomorphic symplectic manifold is torsion-free hence there exists a unique normal projective variety \(\widetilde{Y}_v(\A, \theta)\) with a finite \(2\colon 1\) morphism \(\widetilde{\varepsilon}_v \colon \widetilde{Y}_v(\A,\theta)\to \widetilde{\K}_v(\A,\theta)\) ramified on \(\widetilde{\Sigma}_v\). Moreover there exists a unique normal projective variety \(Y_v(\A,\theta)\) equipped with a finite \(2 \colon 1\) morphism \(\varepsilon_v \colon Y_v(\A, \theta) \to \K_v(\A,\theta)\) whose branch locus is \(\Sigma_v\) \cite[Theorem 4.2]{MRSHodge}. The finite morphism \(\varepsilon_v\) induces a regular involution \(i\) on \(Y_v(\A,\theta)\) hence the morphism \(\varepsilon_v\) can be identified with the quotient map of the involution \(i\). In \cite[Proposition 5.3]{MRSHodge} it is shown that \(Y_v(\A,\theta)\) is always birational to an irreducible holomorphic symplectic manifold of K3\(^{[3]}\) type and a resolution of the indeterminacy of the birational map is explicitly described. We will recall the construction omitting the dependence to the Mukai vector \(v\) and to the abelian surface \(\A\), to avoid cumbersome notation. We denote by \(\Gamma\) the singular locus of \(Y\) which consists of \(256\) points, and by \(\overline{\Gamma}\) the exceptional divisor of \(Bl_{\Gamma}(Y)\) which concerns in the disjoint union of \(256\) copies of the incidence variety; every incidence variety is denoted by \(I_{i}\), and \(I_i \subset \mathbb{P}(V) \times  \mathbb{P}(V)^{\vee}\). Here \(V\) is a 4 dimensional vector space, as we can find in \cite[Section 2]{MRSHodge}. The incidence variety \(I_i \subset \mathbb{P}(V) \times  \mathbb{P}(V)^{\vee}\) has two natural \(\mathbb{P}^{2}\) fibrations given by the projections onto \(\mathbb{P}(V)\) and \(\mathbb{P}(V)^{\vee}\). Let \(p_{i}:I_{i} \longrightarrow \mathbb{P}(V)\) be the two projections. We know that 
the normal bundle of \(I_i\) in \(Bl_{\Gamma}Y\) has degree -1 on the fibers of \(p_i\). Using Nakano's contraction Theorem, \cite{nakano1971inverse}, there exists a complex manifold \(\underline{Y}\) and a morphism of complex manifolds \(h \colon Bl_{\Gamma}Y \longrightarrow \underline{Y}\), whose exceptional locus is \(\overline{\Gamma}\) and such that the image \(J_{i}=h(I_i)\) of any component of \(\overline{\Gamma}\) is isomorphic to \(\mathbb{P}^{3}\). If we consider the restriction of \(h\) to \(I_i\), this is equal to \(p_i\), and \(h\) realizes \(Bl_{\Gamma}Y\) as the blow up of \(\underline{Y}\) along the disjoint union \(J=h(\overline{\Gamma})\) of all the \(J_i\). Moreover by \cite[Proposition 5.3]{MRSHodge} the manifold \(\underline{Y}\) is an irreducible holomorphic symplectic manifold of K3\(^{[3]}\) type. By construction \(\underline{Y}\) has a natural and regular morphism to \(Y\) that contracts \(J\) to \(\Gamma\). Moreover the regular involution \(i\) on \(Y\) is lifted to a birational symplectic involution \(\underline{i}\) on \(\underline{Y}\) which can not be extended to a regular involution \cite[Remark 5.4]{MRSHodge}. More precisely the birational symplectic involution \(\underline{i}\) on \(\underline{Y}\) is regular on the complement of the 256 \(\mathbb{P}^3\)'s.
 In the following \(\widetilde{\K}\) is the manifold of \(\OG_6\) type obtained as resolution of singularities of \(\K\), which is a singular moduli space of sheaves on \(\A\), \(Y\) is the normal projective variety which is singular in 256 points, and \(\underline{Y}\) is the irreducible holomorphic symplectic manifold of K3\(^{[3]}\) type birational to \(Y\). 
Following the original notation of \cite{MRSHodge} we call \(\Delta \subset Y\) the ramification locus (with the reduced induced structure) of \(\varepsilon \colon Y \to \K\). The double cover \(\varepsilon\) induces an isomorphism \(\Delta \cong \Sigma\) between the ramification locus and the singular locus \(\Sigma\) of \(\K\). There exists the following commutative diagram.

\begin{equation}\label{Diagram}
	\begin{tikzcd}
		& & \underline{Y'} \arrow[dd, dashed, "m", leftrightarrow] \\
		& Bl_{\Gamma}Y\arrow[loop left, "\overline{i}" ]\arrow[dd] \arrow[dr, "h_1"] \arrow[ur, "h_2"] &\\
		& & \underline{Y}   \arrow[loop right, dashed, "\underline{i}"]   \arrow[dl]      \\
		& Y \arrow[loop left, "i" ]  \arrow[dd, "\varepsilon" left, "2:1" right] 
		\arrow[ur, bend right, dashed] & \\		\widetilde{\K} \arrow[dr] & & \\
		& \K \stackrel{\text{bir}}{\cong} Y/i& 
	\end{tikzcd}
\end{equation}

\begin{remark}
The codimension of the family of manifolds of \(\OG_6\) type that are moduli space of sheaves on an abelian surface is 3 in the moduli space of marked manifolds of \(\OG_6\) type, denoted by \(\mathcal{M}_{\OG_{6}}\). In fact the Néron–Severi group of a generic element in this family is at least three dimensional, since it contains the class of the exceptional divisor, the class of the locus of non-locally free sheaves and the class arising from the ample divisor \(\theta\) on the abelian surface \(\A\).
A natural question is what is the dimension of the family of manifolds of \(\OG_6\) type that admit a Mongardi--Rapagnetta--Sacc\'{a} model.
The moduli space  \(\mathcal{M}_{K{3}^{[3]}}\) is the marked moduli space of manifolds of \(K{3}^{[3]}\) type which has dimension \(21=h^{1,1}(K3^{[3]})\). 
The manifold \(\underline{Y}\) is a manifold of K3\(^{[3]}\) type and \(i\) is a birational involution defined on it. This involution \(i\) is symplectic hence if \(\sigma_{\underline{Y}}\) is the symplectic form then \(i^{*}\sigma_{\underline{Y}} = \sigma_{\underline{Y}}\) which means that \(\sigma_{\underline{Y}} \in \mathbb{P}( \HH^{2}(\underline{Y}, \mathbb{C})^{i})\) which is a six dimensional complex space.
Since \(\sigma_{\underline{Y}}\) is a symplectic form, \(\sigma_{\underline{Y}} \overline{\sigma_{\underline{Y}}} = 0\), hence it verifies a quadratic equation in a space of dimension six, which means that
\[\{X \mbox{ of \(\OG_6\) type that admit a MRS model}\} \subseteq \mathcal{M}_{\OG_{6}}\]
is a five dimensional subspace of the six-dimensional marked moduli space of \(\OG_6\) type manifolds. Due to this fact it would be possible to generalize the MRS construction for manifolds of \(\OG_6\) type in a codimension 1 subspace of \(\mathcal{M}_{\OG_{6}}\).
\end{remark}

\subsection{Proof of \autoref{ind at the quot} and more remarks}\label{proof thm 1.3}

We introduce the notion of automorphisms or birational transformations induced at the quotient in order to find a criterion to determine when an automorphism or a birational transformation of \(\widetilde{\K}\) lifts to a birational transformation of the manifold of K3\(^{[3]}\) type involved in the construction recalled in \autoref{The Mongardi--Rapagnetta--Sacca model}.
\begin{definition}\label{def aut ind at the quotient}
If \(\widetilde{\K}\) is an irreducible holomorphic symplectic manifold of \(\OG_6\) type obtained as a resolution of moduli space of sheaves on an abelian surface, and if \(\varphi \in \Aut(\widetilde{\K})\) (or \(\varphi \in \Bir(\widetilde{\K}))\) is an automorphisms (or a birational transformation) of \(\widetilde{\K}\), then \(\varphi\) is \textit{induced at the quotient} if \(\varphi\) lifts to a birational transformation of \(\underline{Y}\), where \(\underline{Y}\) is the smooth irreducible holomorphic symplectic manifold of K3\(^{[3]}\) type of diagram \ref{Diagram}.   
\end{definition}

\begin{proposition}\label{teorema sul sollevamento alla birazionalità della K3[3]}
If \(\widetilde{\K}\) is a manifold of \(\OG_6\) type as in diagram \ref{Diagram}, and if \(\varphi \in \Bir(\widetilde{\K})\) is a birational transformation of finite order of \(\widetilde{\K}\) such that there exists a class \(E \in \NS(\widetilde{\K})\) of \((1,1)\)-type, of square \(-2\) and divisibility \(2\) which is fixed by the induced action of \(\varphi\) in cohomology, then \(\varphi\) is induced at the quotient.
\end{proposition}

\proof
Since \(E \in \NS(\widetilde{\K})\) is fixed by the induced action of \(\varphi\) and since the class \(E\) represent the cohomology class of the exceptional divisor of the resolution \(\widetilde{\K} \to \K\), then the automorphism \(\varphi\) is well defined on \(\K\). From \cite[Remark 3.2, Theorem 4.2]{MRSHodge}, we have that if \(\varepsilon \colon Y \longrightarrow \K\) is the étale double cover, then \(\varepsilon^{-1}(\K \ \backslash \ \Sigma) = Y \ \backslash \ \Delta\). Since the real codimension of \(\Delta\) is \(2\), see \cite{MRSHodge}, then the map \(\pi_1(Y\ \backslash \ \Delta) \twoheadrightarrow \pi_1(Y)\) is surjective. We have that \(\pi_{1}(Y \ \backslash \ \Delta) = 0\) and \(\varepsilon : Y \ \backslash \ \Delta \longrightarrow \K \ \backslash \ \Sigma\)  is an étale cover.
We can consider the following diagram: 

\begin{center}	
	\begin{tikzcd}
		Y \ \backslash \ \Delta  \arrow[r, "\widetilde{\psi}"] \arrow[d, "\varepsilon" left, "2:1" right]
		& Y \ \backslash \ \Delta \arrow[d, "\varepsilon" left, "2:1" right] \\
		K \ \backslash \ \Sigma \arrow[r, "\varphi"]
		&K \ \backslash \ \Sigma 
	\end{tikzcd}
	
\end{center}

By \cite[Proposition 1.33]{hatcher2005algebraic} we know that if 
\begin{equation}\label{equation lemma 4.10}
\varphi (\varepsilon ( \pi_{1}(Y \ \backslash \ \Delta) ) )\subseteq \varepsilon(\pi_{1}(Y \ \backslash \ \Delta))
\end{equation}
then \(\varphi\) lifts to an automorphism \(\psi \colon Y \  \backslash \ \Delta \longrightarrow Y \  \backslash  \ \Delta\).  In our case \(\pi_{1}(Y \setminus \Delta) = 0\) hence \eqref{equation lemma 4.10} is verified. The set \(Y \setminus \Delta\) is an open subset of \(Y\) hence \(\psi\) is a birational transformation of \(Y\). Moreover the manifold \(Y\) is birational to the irreducible holomorphic symplectic manifold of \(K3^{[3]}\) type \(\underline{Y}\) hence \(\varphi\) is induced at the quotient.
\endproof

\proof (of \autoref{ind at the quot})
If \(X\) is a numerical moduli space space, then by \autoref{X n.m.s allora birazionale alla fibra della mappa di albanese} there exists an abelian surface \(\A\), a Mukai vector \(v\), and a polarization \(\theta\) on \(\A\), such that there exists a birational map \(\alpha: X \dashrightarrow \widetilde{\K}_v(\A, \theta)\), where \(\widetilde{\K}_v(\A, \theta)\) denotes the resolution of the fiber of the moduli space of sheaves on the abelian surface \(\A\). We denote \(\widetilde{\K}_v(\A, \theta)\) shortly by \(\widetilde{\K}\). Since \(\varphi \in \Bir(X)\) is a birational transformation of \(X\), then \(\widetilde{\varphi}=\alpha \circ \varphi \circ \alpha^{-1} \in \Bir(\widetilde{\K})\) is a birational transformation of the moduli space \(\widetilde{\K}\). By assumption there exists a class \(E \in \NS(X)\) of square \(-2\) and divisibility \(2\) which is preserved by the action of \(\varphi\) hence the same holds true for \(\widetilde{\varphi}\). By  \autoref{teorema sul sollevamento alla birazionalità della K3[3]} we get the result.
\endproof

In the next proposition we prove that actually we can extend \(\psi\colon Y \ \backslash \ \Delta \longrightarrow Y \ \backslash \ \Delta\) to an automorphism of \(Y\). 

\begin{lemma}\label{lemma2}
	In the previous notations, let \(\varepsilon\colon Y \longrightarrow \K\) be the \(2\colon 1\) cover described above, and let \(\varphi\in \Aut(\K)\) be an automorphism of \(\K\). Suppose there exists an open subset \(U\) of \(\K\) such that \(\varphi_{|_U}:U \rightarrow U\) lifts to \(\psi\colon \varepsilon^{-1}(U) \rightarrow \varepsilon^{-1}(U)\), then \(\psi\) extends to a regular morphism \(\widetilde{\psi}:\varepsilon^{-1}(\K) \longrightarrow \varepsilon^{-1}(\K)\) such that \(\widetilde{\psi}_{|_{\varepsilon^{-1}(U)}}=\psi\).
	\begin{center}	
		\begin{tikzcd}
			\varepsilon^{-1}(U) \subseteq Y \arrow[r, "\psi"] \arrow[d, "\varepsilon" left, "2:1" right]
			& \varepsilon^{-1}(U) \subseteq Y \arrow[d, "\varepsilon" left, "2:1" right] \\
			U \subseteq K  \arrow[r, "\varphi"]
			&U \subseteq K 
		\end{tikzcd}
	\end{center}
\end{lemma}

\proof

From hypothesis we know that $\varphi\colon\K \longrightarrow \K$ is regular. If we denote $\Gamma_{\varphi} \subset \K \times \K$ the graph of the morphism, then it is well known that $p_1\colon\Gamma_{\varphi}\xrightarrow{\cong}\K$ is an isomorphism. For the same reason we have the graph $$\Gamma_{\psi} \subset \varepsilon^{-1}(U) \times \varepsilon^{-1}(U)$$ and the isomorphism $p_{1}:\Gamma_{\psi} \xrightarrow{\cong} \varepsilon^{-1}(U) $. We have that $$\Gamma_{\psi} \subset \varepsilon^{-1}(U) \times \varepsilon^{-1}(U) \subseteq Y \times Y $$ where the last is an inclusion in a compact set. We can consider the Zariski closure of the graph, that we denote with $\overline{\Gamma_{\psi}}$. The closure $\overline{\Gamma_{\psi}}$ lies in a closed subset of $Y \times Y$, which is the fiber product over $\K$. To be more precise the fiber product is $Y \times_{\varepsilon, \varphi \circ \varepsilon} Y \subset Y \times Y$.
In the following diagram we denote $Y \times_{\varepsilon, \varphi \circ \varepsilon} Y $ with $\overline{Y \times Y}$.

\begin{center}
	\begin{tikzcd}
		\overline{\Gamma_{\psi}} \arrow[r, hook] \arrow[dr, "\xi"] 
		\arrow[d , dashrightarrow]
		& \overline{Y \times Y}   \arrow[d, "\simeq"] \arrow[r, hook] & Y \times Y \arrow[d, "\overline{\varepsilon}"]\\
		Y \arrow[r, "\varepsilon" ]&\Gamma_{\varphi} \cong K \arrow[r, hook] & K \times K  \end{tikzcd}
\end{center}

In this commutative diagram $\overline{\varepsilon}$ is generically finite, $\overline{\Gamma_{\psi}}$ is a subset of $Y \times_{\varepsilon, \varphi \circ \varepsilon} Y $ and$$Y \times_{\varepsilon, \varphi \circ \varepsilon} Y  \xrightarrow{\cong} \Gamma_{\varphi}$$ is an isomorphism by construction. For this reason $\xi: \overline{\Gamma_{\psi}}  \longrightarrow \K$ is a finite morphism and consequently $ \overline{\Gamma_{\psi}} \longrightarrow Y$ is a finite morphism. Now by hypothesis we have that the previous map is injective on an open subset. Since $Y$ is a normal variety we can conclude that $\overline{\Gamma_{\psi}} \longrightarrow Y$ is an isomorphism, which implies that $\widetilde{\psi} \colon Y \longrightarrow Y$ is a regular morphism, where \(\widetilde{\psi}\) is such that \(\widetilde{\psi}_{|_{\varepsilon^{-1}(U)}}=\psi\).
\endproof
\begin{proposition}\label{induced an aut of Y}
If \(\widetilde{\K}\) is a manifold of \(\OG_6\) type as in diagram \ref{Diagram}, and if \(\varphi \in \Aut(\widetilde{\K})\) is an automorphism of finite order of \(\widetilde{\K}\) such that there exists a class \(E\) of \((1,1)\)-type, of square \(-2\) and divisibility \(2\) which is fixed by the induced action of \(\varphi\) in cohomology, then \(\varphi\) lifts to an automorphism of \(Y\).
\end{proposition}
\proof
Consider \(\varphi \in \Aut(\widetilde{\K})\) an automorphism of the O'Grady's sixfold \(\widetilde{\K}\) then by \autoref{teorema sul sollevamento alla birazionalità della K3[3]} we know that \(\varphi\) lifts to a birational transformation of \(Y\). Moreover taking \(U=\K \setminus \Sigma\) by \autoref{lemma2} we know that \(\varphi\) lifts to an automorphism of \(Y\).
\endproof
In the following diagram we denote with the usual notation the automorphisms and the varieties involved in the construction.

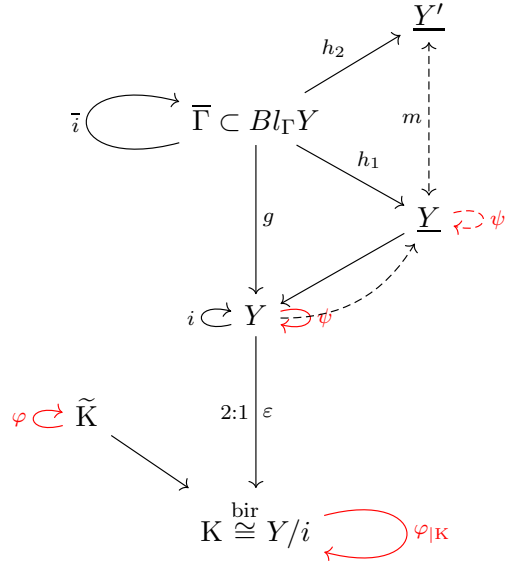
\begin{center}
	\begin{tikzcd}\label{Diagram_maps}
		& & \underline{Y'} \arrow[dd, dashed, "m" left, leftrightarrow] \\
		& \overline{\Gamma} \subset Bl_{\Gamma}Y\arrow[loop left, "\overline{i}" ]     \arrow[dd, "g"] \arrow[dr, "h_1"] \arrow[ur, "h_2"] &\\
		& & \underline{Y}   \arrow[loop right, dashed, red,  "\psi"]   \arrow[dl]      \\
		& Y \arrow[loop left, "i" ]  \arrow[loop right, red , "\psi"]\arrow[dd, "2:1"  left , "\varepsilon", right] 
		\arrow[ur, bend right, dashed] & \\
		\widetilde{\K} \arrow[loop left, red, "\varphi"] \arrow[dr] & & \\
		& \K \stackrel{\text{bir}}{\cong} Y/i \arrow[loop right, red, "\varphi_{|\K}"]& 
	\end{tikzcd}
\end{center}
\subsection{A sufficient condition to have a regular morphism on the Hilbert scheme}\label{A sufficient cond to have a regular morphism of K3}
Now we want to find a criterion to determine when a birational transformation of \(\underline{Y}\) extends to an automorphism of \(\underline{Y}\). Equivalently we want to determine when an automorphism \(\psi\) of the singular normal projective variety \(Y\) lifts to an automorphism \(\underline{\psi}\) of the smooth irreducible holomorphic symplectic manifold \(\underline{Y}\) of K3\(^{[3]}\) type.
In diagram at the end of \autoref{proof thm 1.3} \(\Gamma\) is the singular locus of \(Y\) and it consists of 256 points. We have that \(\psi(\Gamma) =\Gamma\). The automorphism \(\psi\) of \(Y\) preserves the singular locus but can permute the singular points. If we assume that the 256 singular points of \(\Gamma\) are pointwise fixed then the automorphism \(\psi:Y \longrightarrow Y\) extends in a direct way on the blow up of these singular points, which means that \(\overline{\psi}: Bl_{\Gamma}Y \longrightarrow  Bl_{\Gamma}Y\) is a well defined automorphism. What we need to find is a sufficient condition to extend this automorphism on \(\underline{Y}\). As we know by \cite{MRSHodge}, the preimage \(g^{-1}(\Gamma)=\overline{\Gamma}\) is the exceptional divisor of \(Bl_{\Gamma}(Y)\), and consists of 256 copies of the incidence variety \(I_{i}\).

The exceptional locus of \(h_1 \colon Bl_{\Gamma}Y \longrightarrow \underline{Y}\) is \(\overline{\Gamma}\) and the image \(J_{i}=h_1(I_i)\) of any component of \(\overline{\Gamma}\) is isomorphic to \(\mathbb{P}^{3}\).
The automorphism \(\overline{\psi}\) on \(Bl_{\Gamma}(Y)\) descends to a birational map that is well defined outside \(J\), i.e. outside the disjoint union of 256 copies of \(\mathbb{P}^{3}\).
We want to find sufficient conditions to extend this map on these \(\mathbb{P}^{3}\)'s and to obtain an automorphism of \(\underline{Y}\). The preimage with respect to \(g \colon Bl_{\Gamma}Y \to Y\) of a a singular point \(p\) is an incidence variety \(I_i\), which is a divisor of \(Bl_{\Gamma}Y\). By \cite{MRSHodge} we know that on every incidence variety is defined a fibration with basis \(\mathbb{P}^{3}\) and fiber isomorphic to \(\mathbb{P}^{2}\) and there exists the following diagram, see \cite{MRSHodge}.
%
\begin{center}
	\begin{tikzcd}[row sep=huge]
		& \mathbb{P}^{3} \subset \underline{Y'} \arrow[dd, dashed, "Mukai flop" , leftrightarrow]  \\
		\mathbb{P}^3 \times \mathbb{P}^3 \supset I_1 \simeq I_2 \subset Bl_{\Gamma}Y \arrow[ur, "p_{1}" description]  \arrow[dr, "p_{2}" description] &     \\
		& \mathbb{P}^{3} \subset \underline{Y}
	\end{tikzcd}
\end{center}
We call the incidence variety \(I\). Since \(Bl_{\Gamma}Y \cong  Bl_{J}\underline{Y}\), we have the following result.
\begin{proposition}
	The incidence variety \(I\) is isomorphic to \(\mathbb{P}(\Omega_{\mathbb{P}^{3}})\) and \[\Pic(I) \cong \Pic(\mathbb{P}^3 \times (\mathbb{P}^3)^{\vee})  \cong \langle H_1, H_2 \rangle\] where \(H_1 = p_1^{*}( \mathcal{O}_{\mathbb{P}^{3}}(1))\) and \(H_2 = p_2^{*}( \mathcal{O}_{(\mathbb{P}^{3})^{\vee}}(1))\). 
\end{proposition}

\proof
The variety $\underline{Y}$ is an irreducible holomorphic symplectic manifold of dimension 6 and $\mathbb{P}^{3}$ is a lagrangian subspace of $\underline{Y}$. The symplectic form $\sigma_{\underline{Y}}$ gives a duality between $\mathcal{T}_{\mathbb{P}^3} $ and $\Omega_{\mathbb{P}^3}$, but $\sigma_{\underline{Y}}$ on the tangent bundle is zero; this duality is the one that sends $\mathcal{N}_{\mathbb{P}^3}$ to $\Omega_{\mathbb{P}^3}$ which are isomorphic. We know that the exceptional locus of this blow up is $I \cong \mathbb{P}(\mathcal{N}_{\mathbb{P}^3}) \cong \mathbb{P}(\Omega_{\mathbb{P}^3})$. We define $\mathcal{O}_{\mathbb{P}^{3}}(1) \boxtimes \mathcal{O}_{\mathbb{P}^{3}}(1) := p_1^{*}( \mathcal{O}_{\mathbb{P}^{3}}(1)) \otimes p_2^{*}( \mathcal{O}_{\mathbb{P}^{3}}(1))$.
Since on $I$ two $\mathbb{P}^{2}$ fibrations are well defined, if we call $H_1=p_1^{*}( \mathcal{O}_{\mathbb{P}^{3}}(1))$ and $H_2=p_2^{*}( \mathcal{O}_{\mathbb{P}^{3}}(21))$ we can say that $ \Pic(\mathbb{P}^{3} \times \mathbb{P}^{3})$ is generated by $\mathcal{O}_{\mathbb{P}^{3}}(1) \boxtimes \mathcal{O}_{\mathbb{P}^{3}}(1)$. By Lefschetz's  Theorem for the Picard group, we know that $\Pic(I)=\Pic(\mathbb{P}^{3} \times \mathbb{P}^{3})=\langle H_1, H_2 \rangle$, where $H_{1}$ comes from the first fibration and $H_{2}$ comes from the second fibration.
\endproof
In the next theorem we show that a sufficient condition for an automorphism \(\overline{\psi}\) of \(Bl_{\Gamma}Y\) to descend to an automorphism of \(\underline{Y}\) is to not exchange the fibers of the two \(\mathbb{P}^{2}\) fibrations. 


\begin{theorem}\label{thm induced at the quotient}
	Let \(\widetilde{\K}\) be a manifold of \(\OG_6\) type obtained as a resolution of a moduli space. Let \(\varphi \in \Aut(\widetilde{\K})\) be an automorphism of prime order \(p\), \(p \neq 2\), and suppose there exists a class \(E \in \NS(\widetilde{\K})\) with \(E^{2}=-2\) and divisibility \(2\) which is preserved by the induced action of \(\varphi\), then \(\varphi\) is induced at the quotient. Denote by \(\psi\) the lifted action on \(Y\) and assume that the 256 singular points of \(Y\) are pointwise fixed by \(\psi\), then \(\psi\) lifts to an automorphism \(\underline{\psi}\) of \(\underline{Y}\).
\end{theorem}
To give a proof of \autoref{thm induced at the quotient} we need the following lemma.
\begin{lemma}\label{lemma Pic(I)=H_1,H_2}
	Let $f$ be an automorphism of $Bl_{\Gamma}Y$ 
	that leaves invariant the exceptional divisor,
	and $f^{\ast}$ the induced action on $\Pic(I)=\langle H_1, H_2 \rangle $. Then $f^{*}$ is the identity or $f^{*}(H_1)=H_2$ and $f^{*}(H_2)=H_1$.
\end{lemma}
\proof
In the proof we denote \(\mathcal{O}_{\mathbb{P}^3}(1)\) by \(H\) and \(\mathcal{O}_{(\mathbb{P}^3)^{\vee}}(1)\) by \(H^{\vee}\). 
By assumption \(H_{1}=p_1^{*}(H)\) and \(H_{2}=p_2^{*}(H^{\vee})\) are the generators of the hyperplane sections of \(p_1\) and \(p_2\). Note that \(H_1\) corresponds to the cycle \([(H \times (\mathbb{P}^3)^{\vee})\cap I]\) and \(H_2\) corresponds to the cycle \([(\mathbb{P}^3 \times H^{\vee})\cap I]\) where the class is in the Chow group.
Moreover $H_1^2$ is the class corresponding to the cycle $[(l     \times (\mathbb{P}^3)^*)\cap I]$, 
where the class is in the Chow group. Moreover, for $H_2$ it holds the same: $H_2^3$ is the class corresponding to the cycle $[( \mathbb{P}^3 \times  p)\cap I]$.
This is the fiber of the closed point $p$
and this is isomorphic to $\mathbb{P}^{2}$.
The product $H_1^{2}H_2^{3}$ is equal to 1, since this is an intersection of a line and a $\mathbb{P}^{2}$ in a generic position.
With the same argument, but exchanging the role of $H_1$ and $H_2$ we obtain that $H_1^{3}H_2^{2}$ is equal to 1. The pullback commutes with the intersection product, hence for dimensional reasons the product \(H_1^{k}\) is equal to zero when \(k \geq 4\) and the same holds true for \(H_2\).
Moreover, since the pullback operation commutes with the intersection form,
we have that $f^*(H_1)^5=f^{*}(H_1^5)=0$. The action of $f^{*}$ preserves the Picard group of $I$, hence we can denote $f^{*}(H_1)=\alpha H_1 + \beta H_2$ and $f^{*}(H_2)=\gamma H_1 + \delta H_2$. With this notation we have: 
$$(f^{*}H_1)^5=\sum_{i=0}^5 \binom{5}{i}  \alpha^i \beta^{5-i} H_1^{i}H_2^{5-i}=10 \alpha^{2}\beta^{3}H_1^{2}H_2^{3}+10 \alpha^{3}\beta^{2}H_1^{3}H_2^{2}= 10 \alpha^{2}\beta^{3}+10 \alpha^{3}\beta^{2
}.$$
Furthermore we have $$\alpha^2 \beta^2 (\alpha + \beta)=0.$$
In the same way for $H_2$ we obtain: $$\gamma^2 \delta^2 (\gamma + \delta)=0.$$
After some straightforward computation we obtain the following six cases:
\begin{multicols}{3}
	$\begin{cases}
	f^{*}(H_1)=H_1 \\ f^{*}(H_2)=H_2
	\end{cases}$

	$\begin{cases}
	f^{*}(H_1)=\pm(H_1-H_2)\\ f^{*}(H_2)=H_2
	\end{cases}$
	
	$\begin{cases}
	f^{*}(H_1)=H_1 \\ f^{*}(H_2)=\pm(H_1-H_2)
	\end{cases}$
	
\end{multicols}

\begin{multicols}{3}
	
	$\begin{cases}
	f^{*}(H_1)=H_2 \\ f^{*}(H_2)=\pm(H_1-H_2)
	\end{cases}$

	$\begin{cases}
	f^{*}(H_1)=H_2 \\ f^{*}(H_2)=H_1
	\end{cases}$

	$\begin{cases}
	f^{*}(H_1)=\pm(H_1-H_2)\\ f^{*}(H_2)=H_1
	\end{cases}$

\end{multicols}
We can notice that $f^{*}(H_1)=\pm(H_1-H_2)$ is not allowed. In fact, let $l_{1} \subset p_{1}^{-1}(p) \simeq \mathbb{P}^{2}$ and let $l_{2} \subset p_{2}^{-1}(q) \simeq \mathbb{P}^{2}$ be two lines which lie in the two different fibrations. Since $f^{\ast}H_1 . l_1= f_{\ast}(f^{\ast}H_{1}. l_1)=H_{1}. f_{\ast}l_1$, we notice that $f_{\ast} l_1$ is a line which means $f_{\ast} l_1 \cong \mathbb{P}^{1}$ since $f$ is an automorphism and for this reason $H_{1}. f_{\ast}l_1$ could be 1 or 0. Assume that $H_{1}. f_{\ast}l_1=1$ and that $f^{*}(H_1)=H_1-H_2$. We want to compute the intersection \((H_1-H_2)\cdot l_1=H_1 \cdot l_1 - H_2 \cdot l_1\). It holds that \(H_1=[(H \times (\mathbb{P}^3)^{\vee})\cap I]\) and \(l_1=[(p \times l_1^{\vee}) \cap I]\) where \(l_1=p_1^{-1}(p)\). The point \(p\) and the hyperplane \( H\) are generically disjoint in \(\mathbb{P}^3\), hence \(H_1 \cdot l_1=0\). On the other hand \(H_2=[(\mathbb{P}^3 \times H^{\vee})\cap I]\) hence \(H_2 \cdot l_1=1\) by the generic position of the hyperplane \(H^{\vee}\) and the line \(l_1^{\vee}\) in \((\mathbb{P}^3)^{\vee}\). We deduce that \((H_1-H_2)\cdot l_1=H_1 \cdot l_1 - H_2 \cdot l_1=-1\) that is an absurd.
This holds true also in the other similar cases using properly \(l_1\) or \(l_2\) and \(H_1\) or \(H_2\). We can conclude that the two possible actions of $f^{*}$ on $\Pic(I)$ are the identity and the automorphism that exchanges $H_1$ and $H_2$.
\endproof

\proof (of \autoref{thm induced at the quotient}) By \autoref{induced an aut of Y} we know that \(\varphi \in \Aut(\widetilde{\K})\) lifts to an automorphism \(\psi\) on \(Y\) hence \(\varphi\) lifts to a birational transformation of \(\underline{Y}\) and this implies that \(\varphi\) is induced at the quotient. By assumption the singular points of \(Y\) are pointwise fixed hence \(\psi\) lifts to \(\overline{\psi}\). By \autoref{lemma Pic(I)=H_1,H_2} we know the action of \(\overline{\psi}\) on \(\Pic(I)\), hence we deduce that if the order of the automorphism is prime \(p\), with \(p > 2\), the action is the identity on \(\Pic(I)\). As a consequence the fibers of the two \(\mathbb{P}^{2}\) fibrations are not exchanged and we can define an automorphism \(\underline{\psi}\) on \(\underline{Y}\).
\endproof

\section{Applications}\label{Applications}

In this section we give an application of the two criteria for induced automorphisms described in \autoref{Induced automorphisms groups} and \autoref{Automorphisms induced at the quotient}. The lattice-theoretic criterion in \autoref{Induced automorphisms groups} allows to determine if an automorphism or a birational transformation of a manifold of \(\OG_6\) type which is at least birational to \(\widetilde{\K}_v(\A, \theta)\), is induced by an automorphism of the abelian surface \(\A\). 

Similarly to determine if an automorphism or a birational transformation of a manifold \(\widetilde{\K}_v(\A,\theta)\) lifts to a birational transformation of the manifold of K3\(^{[3]}\) type involved in the MRS construction it is enough to know its induced action on the second integral cohomology lattice. 

\begin{proposition}\label{induced implies induced at the quotient}
Let \(X\) be a manifold of \(\OG_6\) type which is a numerical moduli space, and let \(\G \in \Bir(X)\) be a finite group of induced birational transformations, then \(\G \subset \Bir(X)\) is a group of birational transformations induced at the quotient.
\end{proposition}
\proof
If \(X\) is a numerical moduli space, then by \autoref{X n.m.s allora birazionale alla fibra della mappa di albanese} \(X\) is birational to the moduli space \(\widetilde{\K}_v(\A, \theta)\). 
If the group \(\G\) is induced then by \autoref{G induced implies G num ind} it is numerically induced, hence there exists a class of \((1,1)\)-type, of square \(-2\) and divisibility \(2\), invariant with respect to the action of \(\G\). By \autoref{ind at the quot} \(\G\) is induced at the quotient.
\endproof

We consider the classification of nonsymplectic automorphisms of prime order on manifolds of \(\OG_6\) type given in \cite[Table 1]{Grossi:non-sp.aut.OG3}, where the author considers a manifold of \(\OG_6\) type, a fixed marking \(\eta:\HH^{2}(X,\mathbb{Z}) \to \bL\) of \(X\), and classifies the invariant and the coinvariant sublattices, denoted by \(\bL^{\G}\) and \(\bL_{\G}\) respectively, with respect to the induced action on the second integral cohomology lattice by a nonsymplectic automorphism of prime order. More precisely, the images of nonsymplectic automorphisms of prime order of the representation map 
\[
\eta_{*}:\Aut(X) \to \OO(\bL)
\]
\[ f \mapsto \eta \circ f^{*} \circ \eta^{-1}
\]
are classified.
We use the lattice-theoretic criterion of \autoref{num ind impliesd ind} to determine if a nonsymplectic automorphism of a manifold of \(\OG_6\) type is induced. Similarly we use the lattice-theoretic criterion of  \autoref{ind at the quot} to determine if an automorphism is induced at the quotient.

\begin{remark}\label{remark sul NS}
If \(\varphi\in \Aut(X)\) is a nonsymplectic automorphism of prime order \(p\), by \cite[Proposition 3.3]{Grossi:non-sp.aut.OG3} we know that \(p \in \{2,3,5,7\}\). Moreover if we denote by \(\G\) the cyclic group of prime order generated by \(\varphi\), then by \cite[Remark 3.2]{Grossi:non-sp.aut.OG3} we can assume that \(\bL^{\G}=\NS(X)\) and \(\bL_{\G}=\T(X)\) is the transcendental lattice.
\end{remark}
\begin{remark}
In \autoref{tab:L} there is a classification of invariant and coinvariant sublattices with respect to a nonsymplectic automorphism of prime order on a manifold of \(\OG_6\) type \cite[Table 1]{Grossi:non-sp.aut.OG3}. If \(|\G^{\sharp}|=2\) then \(\G^{\sharp}\) exchanges the two generators of \(\bL^{\sharp}\), hence there are no vectors \(v \in \bL\) such that \(v^{2}=-2\) and \((v, \bL)=2\) and that are fixed by \(\G\). In this case \(X\) is not a numerical moduli space.
\end{remark}

\subsection{Proof of \autoref{theorem applications}}\label{proof of thm 1.4}
\proof
In \autoref{tab:L} we have a classification of \(\bL^{\G}\). We know by \autoref{remark sul NS} that \(\bL^{\G}=\bL^{1,1}\) hence by \autoref{X n.m.s allora birazionale alla fibra della mappa di albanese}, checking the numerical conditions, we determine if \(X\) is a numerical moduli space.

If \(|\G|=2\) by \autoref{cor 2.2 GOV}, we know that for every \(v \in \bL^{\G}\) it holds that \((v, \bL)=1\). In this way we exclude cases \(1,2,3,6,7,11\). Furthermore in cases \(4\) and \(5\), the manifold \(X\) is not a numerical moduli space because by \autoref{relazione segnature} the signature of \(\bLambda_8^{1,1}\) is \((2,0)\) hence it does not contain \(\bU\) as a direct summand. Moreover in cases \(23\) and \(24\) if \(v\) is a vector in \(\bL^{\G}=\bU \oplus \bD_4(-1)\) or in \(\bL^{\G}=\bU(2) \oplus \bD_4(-1)\) with \(v^2=-2\) then \((v, \bL^{\G})=1\). By \autoref{relazione sulle divisibilità} \((v, \bL)=1\) hence in these cases \(X\) is not a numerical moduli space. 
In case \(8\) by \autoref{relazione segnature} the signature of \(\bLambda_8^{1,1}\) is \((2,1)\). In the nonsymplectic setting \(\bL_{\G}=\T(X)\) and \(\bLambda_8^{1,1}\)
are orthogonal complement in the unimodular lattice \(\bLambda_8\) hence by \autoref{embedding in a unim lattice} \(l(\bL_{\G}^{\sharp})=l((\bLambda_8^{1,1})^{\sharp})=3\). The lattice \(\bLambda_8^{1,1}\) is a \(2\)-elementary lattice, its signature is \((2,1)\) and the length of its discriminant group is \(3\), hence it does not contain any copy of \(\bU\) as a direct summand. We can conclude similarly in case \(13\), since the signature of \(\bLambda_8^{1,1}\) is \((2,2)\) and the length of its discriminant group is \(4\). In cases \(9,10,12,14,15,16,17,18,19,20,21,22\) we know the signature, the quadratic form and the length of the discriminant group of the \(2\)-elementary lattice  \(\bLambda_8^{1,1}\), hence we check by \cite[Theorem 1.5.2]{dolgachev1982integral} that there exists a copy of \(\bU\) in \(\bLambda_8^{1,1}\). Moreover we can compute the gluing subgroup since we know that \(\bL^{\sharp}=(\mathbb{Z}/2\mathbb{Z})^{\oplus2}\). Then, by \autoref{Lemma 2.1 GOV} we find that there exists a vector of square \(-2\) and divisibility \(2\) in \(\bL^{\G}=\NS(X)\) hence the manifold \(X\) is a numerical moduli space, so by \autoref{ind at the quot} the automorphisms are induced at the quotient. 
Between automorphisms that are induced at the quotient we can detect which ones are induced by \autoref{num ind impliesd ind} only checking that the third condition of \autoref{numerically induced} is verified (the other two are already verified). Since we are taking antisymplectic involutions we only need to check that the rank of \(\bL_{\G}\) is even and we deduce that in cases \(12,14,15,16,17,18\) the automorphisms are also induced. 

If \(|\G|=3\) in case \(1\) the manifold \(X\) is not a numerical moduli space because by \autoref{relazione segnature} the signature of \(\bLambda_8^{1,1}\) is \((2,0)\) hence we can not find any copy of \(\bU \) in it as a direct summand. In case \(2\) of \(|\G|=3\) and for \(|\G|=5\) there is a copy of \(\bU\) in \(\bL^{\G}\) hence in \(\bLambda_8^{1,1}\) which is absurd.
Moreover we can compute the gluing subgroup since we know that \(\bL^{\sharp}=(\mathbb{Z}/2\mathbb{Z})^{\oplus2}\). Then, by \autoref{Lemma 2.1 GOV} we find that there exists a vector of square \(-2\) and divisibility \(2\) in \(\bL^{\G}=\NS(X)\), hence \(X\) is a numerical moduli space and the automorphism is induced at the quotient. Moreover in the latter two cases the third condition of \autoref{numerically induced} is verified hence the automorphisms are also induced.

If \(|\G|=7\) the manifold \(X\) is not a numerical moduli space because by \autoref{relazione segnature} the signature of \(\bLambda_8^{1,1}\) is \((2,0)\) hence it do not contain any copy of \(\bU \) as direct summand.
\endproof


\bibliographystyle{plain}
\bibliography{Biblio}
\end{document}